\documentclass[onefignum,onetabnum]{siamart190516}

\usepackage[utf8]{inputenc}
\usepackage{lipsum}
\usepackage{amsfonts, amssymb}
\usepackage{graphicx}
\usepackage{epstopdf}
\usepackage{algorithmic}
\ifpdf
  \DeclareGraphicsExtensions{.eps,.pdf,.png,.jpg}
\else
  \DeclareGraphicsExtensions{.eps}
\fi

\newsiamremark{remark}{Remark}
\newsiamremark{hypothesis}{Hypothesis}
\crefname{hypothesis}{Hypothesis}{Hypotheses}
\newsiamthm{claim}{Claim}

\headers{Fast Canonical Polyadic Decomposition}{F. B. Diniz}

\title{A Fast Implementation for the Canonical Polyadic Decomposition\thanks{Submitted to the editors November 7, 2019.
\funding{The author was partially supported by CAPES, Coordenação de Aperfeiçoamento de Pessoal de Nível Superior (Brazil) (Finance code 001).} This work is based on the author's Phd thesis.}}

\author{F. B. Diniz\thanks{Instituto de Matemática, Universidade Federal do Rio de
Janeiro (felipebottega@gmail.com).}}

\usepackage{amsopn}

\ifpdf
\hypersetup{
  pdftitle={A Fast Implementation for the Canonical Polyadic Decomposition},
  pdfauthor={F. B. Diniz}
}
\fi

\begin{document}

\maketitle

\begin{abstract}
A new implementation of the canonical polyadic decomposition (CPD) is presented. It features lower computational complexity and memory usage than the available state of art implementations available. The CPD of tensors is a challenging problem which has been approached in several manners. Alternating least squares algorithms were used for a long time, but they convergence properties are limited. Nonlinear least squares (NLS) algorithms - more precisely, damped Gauss-Newton (dGN) algorithms - are much better in this sense, but they require inverting large Hessians, and for this reason there is just a few implementations using this approach. In this paper, we propose a fast dGN implementation to compute the CPD. In this paper, we make the case to always compress the tensor, and propose a fast damped Gauss-Newton implementation to compute the canonical polyadic decomposition. 
\end{abstract}



\section{Introduction}
In this paper we introduce \emph{Tensor Fox}, a new tensor package that outperform other well-known available tensor packages with respect to the CPD. The main improvements of Tensor Fox with respect to the other packages are: a preprocessing stage, where the multilinear singular value decomposition (MLSVD) is used to compress the tensor; the approximated Hessian has a block structure which can be exploited to perform fast iterations. We will see that Tensor Fox is competitive, being able to handle every tested problem with high speed and accuracy. The other tensor packages considered are: \emph{Tensorlab \cite{tensorlab, tensorlab2} (version 3.0), Tensor Toolbox \cite{tensortoolbox, cp_opt}} (version 3.1) and \emph{TensorLy} \cite{tensorly}. These packages are the state of art in terms of performance and accuracy with respect to the problem of computing the CPD.
	
For decades tensor decompositions have been applied to general multidimensional data with success. Today they excel in several applications, including blind source separation, dimensionality reduction, pattern/image recognition, machine learning and data mining, see for instance \cite{bro, cichoki2009, cichoki2014, savas, smilde, anandkumar}. In the old days, the alternating least squares (ALS) was the working horse algorithm to compute the CPD of a tensor. It is easy to implement, each iteration is cheap to compute and it has low memory cost. However Kruskal, Harshman and others started to notice that this algorithm had its limitations in the presence of \emph{bottlenecks} and \emph{swamps} \cite{kruskal-harshman, burdick}. Several attempts were made in order to improve the performance of ALS algorithms, but none was fully satisfactory. The dGN algorithm is a nonlinear least squares (NLS) algorithm which appears to be a better approach for the CPD. It is less sensitive to bottlenecks and swamps and appears to attain quadratic convergence when close to the solution. The main challenge with the dGN algorithm is the fact that we need to invert a large Hessian at each iteration. Tensor Fox relies on this second approach.
	
The paper is organized as follows. Notations and basic tensor algebra are briefly reviewed in \cref{sec:notations}. The algorithm used to compute the CPD is explained at some detail in section 3. Experiments to validate the claims about the performance of Tensor Fox (TFX) are described in section 4. As an illustration of a 'real world' problem, more tests with Gaussian mixtures are conducted in section 5. The main outcomes of this paper are highlighted in section 6 where we conclude the paper.	

\section{Notations and basic definitions}
\label{sec:notations}

A tensor is a real multidimensional array in $\mathbb{R}^{I_1 \times I_2 \times \ldots \times I_L}$, where $I_1, I_2, \ldots, I_L$ are positive integer numbers. The \emph{order} of a tensor is the number of indexes associated to each entry of the tensor. For instance, a $L$-th order tensor is an element of the space $\mathbb{R}^{I_1 \times I_2 \times \ldots \times I_L}$. Each index is associated to a dimension of the tensor, and each dimension is called a \emph{mode} of the tensor. Tensors of order $1$ are vectors and tensors of order $2$ are matrices. Scalars are denoted by lower case letters. Vectors are denoted by bold lower-case letters, e.g., $\textbf{x}$. Matrices are denoted by bold capital letters, e.g., $\textbf{X}$. Tensors are denoted by calligraphic capital letters, e.g., $\mathcal{T}$. The $i$-th entry of a vector $\textbf{x}$ is denoted by $x_i$, the entry $(i,j)$ of a matrix $\textbf{X}$ is denoted by $x_{ij}$, and the entry $(i,j,k)$ of a third order tensor $\mathcal{T}$ is denoted by $t_{ijk}$. Finally, any kind sequence will be indicated by superscripts. For example, we write $\textbf{X}^{(1)}, \textbf{X}^{(2)}, \ldots$ for a sequence of matrices. We denote by $\ast$ the \emph{Hadamard product} (i.e., the coordinatewise product) between two matrices with same shape.
	
	Through most part of the text we limit the discussion to third order tensors, so instead of using $I_1, I_2, I_3$ to denote the dimensions we will just use $m, n, p$. However, everything can be easily generalized to $L$-th order tensors. 
	
	\begin{definition} The \emph{inner product} of two tensors $\mathcal{T}, \mathcal{S} \in \mathbb{R}^{m \times n \times p}$ is defined as
	$$\langle \mathcal{T}, \mathcal{S} \rangle = \sum_{i=1}^m \sum_{j=1}^n \sum_{k=1}^p t_{ijk} s_{ijk}.$$
	\end{definition}
	
	\begin{definition} The associate (Frobenius) \emph{norm} of $\mathcal{T}$ is given by
	$$\| \mathcal{T} \| = \sqrt{ \langle \mathcal{T}, \mathcal{T} \rangle }.$$
	\end{definition}
	
	We will use the notation $\| \cdot \|$ also for the matrix and vector Frobenius norm. 
	
	\begin{definition} \label{tensprod} The \emph{tensor product} (also called \emph{outer product}) between three vectors $\textbf{x} \in \mathbb{R}^m, \textbf{y} \in \mathbb{R}^n, \textbf{z} \in \mathbb{R}^p$ is the tensor $\textbf{x} \otimes \textbf{y} \otimes \textbf{z} \in \mathbb{R}^{m \times n \times p}$ such that
	$$(\textbf{x} \otimes \textbf{y} \otimes \textbf{z})_{ijk} = x_i y_j z_k.$$
	\end{definition}
	
	We can consider this definition for only two vectors $\textbf{x, y}$. By doing this we have the tensor $\textbf{x} \otimes \textbf{y} \in \mathbb{R}^{m \times n}$ given by
	$$(\textbf{x} \otimes \textbf{y})_{ij} = x_i y_j.$$
	This tensor is just the rank one matrix $\textbf{xy}^T$, which is the outer product between $\textbf{x}$ and $\textbf{y}$.
	
	\begin{definition} A tensor $\mathcal{T} \in \mathbb{R}^{m \times n \times p}$ is said to have \emph{rank} $1$ if there exists vectors $\textbf{x} \in \mathbb{R}^m, \textbf{y} \in \mathbb{R}^n, \textbf{z} \in \mathbb{R}^p$ such that
	$$\mathcal{T} = \textbf{x} \otimes \textbf{y} \otimes \textbf{z}.$$
	\end{definition}
	
	\begin{definition} A tensor $\mathcal{T} \in \mathbb{R}^{m \times n \times p}$ is said to have \emph{rank} $R$ if $\mathcal{T}$ can be written as sum of $R$ rank one tensors and cannot be written as a sum of less than $R$ rank one tensors. In this case we write $rank(\mathcal{T}) = R$.
	\end{definition}
	
	When $\mathcal{T}$ has rank $R$, there exists vectors $\textbf{x}_1, \ldots, \textbf{x}_R \in \mathbb{R}^m$, $\textbf{y}_1, \ldots, \textbf{y}_R \in \mathbb{R}^n$ and $\textbf{z}_1, \ldots, \textbf{z}_R \in \mathbb{R}^p$ such that
	$$\mathcal{T} = \sum_{r=1}^R \textbf{x}_r \otimes \textbf{y}_r \otimes \textbf{z}_r.$$
	By definition we can't write $\mathcal{T}$ with less than $R$ terms.
	
	\begin{definition} A \emph{CPD of rank $R$} for a tensor $\mathcal{T}$ is a decomposition of $\mathcal{T}$ as a sum of $R$ rank one terms.
	\end{definition}
	
	In contrast to the SVD for the matrix case, the CPD is ``unique'' when the rank of a tensor $\mathcal{T}$ of order $\geq 3$ is $R$. The notion of uniqueness here is modulo permutations and scaling. Note that one can arbitrarily permute the rank one terms or scale them as long as their product is the same. These trivial modifications doesn't change the tensor and we consider it to be the unique up to these modifications. 
	
	In practical applications one may not know in advance the rank $R$ of $\mathcal{T}$. In this case we look for a CPD with low rank \cite{comon} such that
	$$\mathcal{T} \approx \sum_{r=1}^R \textbf{x}_r \otimes \textbf{y}_r \otimes \textbf{z}_r.$$  
	
	\begin{definition} The \emph{multilinear multiplication} between the tensor $\mathcal{T} \in \mathbb{R}^{m \times n \times p}$ and the matrices $\textbf{A} \in \mathbb{R}^{m' \times m}, \textbf{B} \in \mathbb{R}^{n' \times n}, \textbf{C} \in \mathbb{R}^{p' \times p}$ is the tensor $\mathcal{T}' \in \mathbb{R}^{m' \times n' \times p'}$ given by
	$$\mathcal{T}'_{i' j' k'} = \sum_{i=1}^m \sum_{j=1}^n \sum_{k=1}^p a_{i' i} b_{j' j} c_{k' k} t_{ijk}.$$
	\end{definition}
	
	To be more succinct we denote $\mathcal{T}' = (\textbf{A}, \textbf{B}, \textbf{C}) \cdot \mathcal{T}$. Observe that in the matrix case we have $(\textbf{A}, \textbf{B}) \cdot \textbf{M} = \textbf{A} \textbf{M B}^T$.  
	
	\begin{definition} Let $\mathcal{T} \in \mathbb{R}^{m \times n \times p}$ be any tensor. A \emph{Tucker decomposition} of $\mathcal{T}$ is any decomposition of the form
	$$\mathcal{T} = (\textbf{A}, \textbf{B}, \textbf{C}) \cdot \mathcal{S}.$$ 	
	\end{definition}
	$\mathcal{S}$ is called the \emph{core tensor} while $\textbf{A}, \textbf{B}, \textbf{C}$ are called the \emph{factors} of the decomposition. Usually they are considered to be orthogonal. One particular (and trivial) Tucker decomposition of $\mathcal{T}$ is $\mathcal{T} = (\textbf{I}_m, \textbf{I}_n, \textbf{I}_p) \cdot \mathcal{T}$, where $\textbf{I}_m$ is the $m \times m$ identity matrix. One way conclude this is just by expanding $\mathcal{T}$ in canonical coordinates by
	$$\mathcal{T} = \sum_{i=1}^m \sum_{j=1}^n \sum_{k=1}^p t_{ijk} \cdot \textbf{e}_i \otimes \textbf{e}_j \otimes \textbf{e}_k,$$
where $\textbf{e}_i, \textbf{e}_j, \textbf{e}_k$ are the canonical basis vectors of their corresponding spaces.

	If $\mathcal{T}$ has rank $R$, let 
	$$\mathcal{T} = \sum_{r=1}^R \textbf{x}_r \otimes \textbf{y}_r \otimes \textbf{z}_r$$
be a CPD of rank $R$. Then a more useful Tucker decomposition is given by $\mathcal{T} = (\textbf{X}, \textbf{Y}, \textbf{Z}) \cdot \mathcal{I}_{R \times R \times R}$, where $\textbf{X} = [ \textbf{x}_1, \ldots, \textbf{x}_R ] \in \mathbb{R}^{m \times R}$, $\textbf{Y} = [ \textbf{y}_1, \ldots, \textbf{y}_R ] \in \mathbb{R}^{n \times R}$, $\textbf{Z} = [ \textbf{z}_1, \ldots, \textbf{z}_R ] \in \mathbb{R}^{p \times R}$ and $\mathcal{I}_{R \times R \times R}$ is the cubic tensor with ones at its diagonal and zeros elsewhere.

\section{Tensor compression}
\label{sec:compression}

Working with ``raw'' data is, usually, not advised because of the typical large data size. A standard approach is to compress the data before starting the actual work, and this is not different in the context of tensors. This is efficient a way of reducing the effects of the curse of dimensionality. We consider tensors of general order only in this section.

    Given a tensor $\mathcal{S} \in \mathbb{R}^{I_1 \times \ldots \times I_L}$, let $\mathcal{S}_{i_\ell = k} \in \mathbb{R}^{I_1 \times \ldots \times I_{\ell-1} \times I_{\ell+1} \times \ldots \times I_L}$ be the sub tensor of $\mathcal{S}$ obtained by fixing the $\ell$-th index of $\mathcal{S}$ with value equal to $k$ and varying all the other indexes. More precisely, $\mathcal{S}_{i_\ell = k} = \mathcal{S}_{: \ldots : k : \ldots :}$, where the value $k$ is at the $\ell$-th index. We call these sub tensors by \emph{hyperslices}. In the case of a third order tensors, these sub-tensors are the slices of the tensor. $\mathcal{S}_{i_1 = k}$ are the horizontal slices, $\mathcal{S}_ {i_2 = k}$ the lateral slices and $\mathcal{S}_ {i_3 = k}$ the frontal slices. The main tool to compress a tensor is the given by the following theorem. 
	
	\begin{theorem}[L. D. Lathauwer, B. D. Moor, J. Vandewalle, 2000] \label{MLSVD}
		For any tensor $\mathcal{T} \in \mathbb{R}^{I_1 \times \ldots \times I_L}$ there are orthogonal matrices $\textbf{U}^{(1)} \in \mathbb{R}^{I_1 \times R_1}, \ldots, \textbf{U}^{(L)} \in \mathbb{R}^{I_L \times R_L}$ and a tensor $\mathcal{S} \in \mathbb{R}^{R_1 \times \ldots \times R_L}$ such that 
		\begin{enumerate}
			\item For all $\ell = 1 \ldots L$, we have $R_\ell \leq I_\ell$. 
			\item $\mathcal{T} = (\textbf{U}^{(1)}, \ldots, \textbf{U}^{(L)}) \cdot \mathcal{S}$.
			\item For all $\ell = 1 \ldots L$, the sub-tensors $\mathcal{S}_{i_\ell = 1}, \ldots, \mathcal{S}_{i_\ell = I_\ell}$ are orthogonal with respect to each other.
			\item For all $\ell = 1 \ldots L$, we have $\|\mathcal{S}_{i_\ell = 1}\| \geq \ldots \geq \|\mathcal{S}_{i_\ell = I_\ell}\|$.
		\end{enumerate}
	\end{theorem}
	
\begin{figure}[htbp]
		\centering
		\includegraphics[scale=.11]{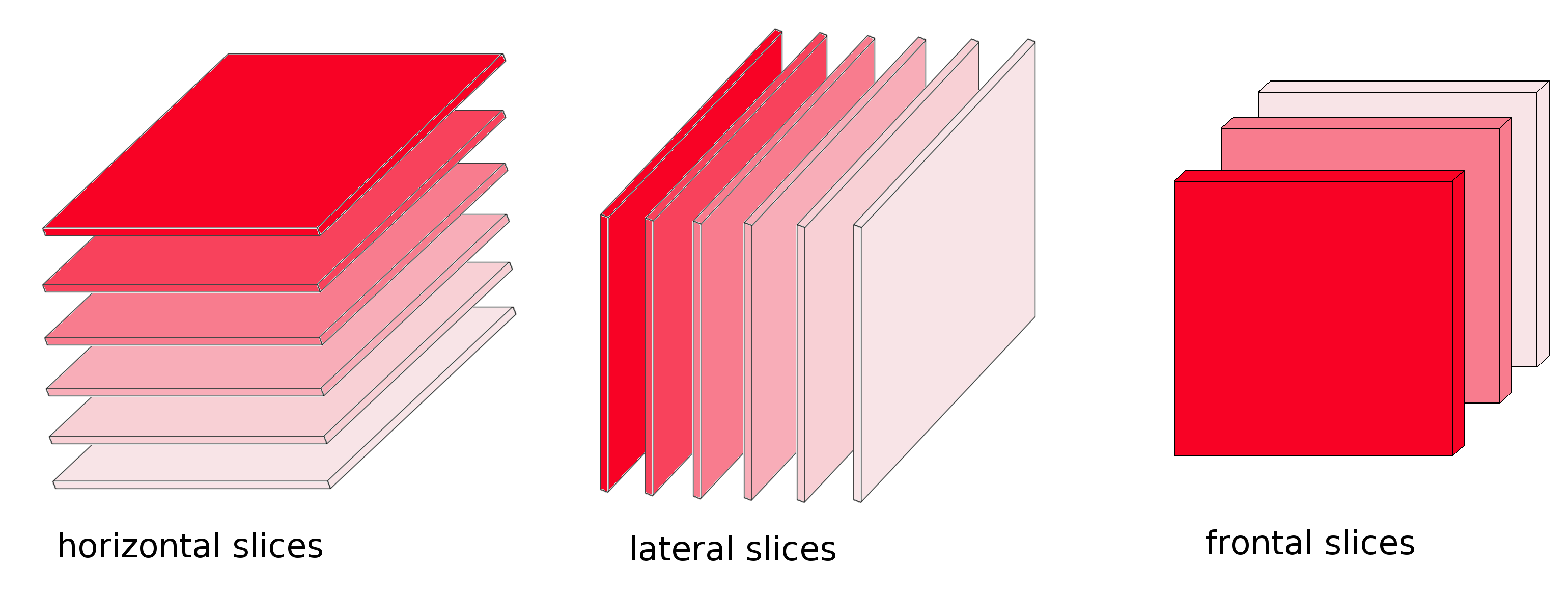}
		\caption{Energy of the slices of a core third order tensor $\mathcal{S}$ obtained after a MLSVD.}
		\label{fig:slices-energy}
	\end{figure}

	We refer to this decomposition as the \emph{multilinear singular value decomposition} (MLSVD) of $\mathcal{T}$. The core tensor $\mathcal{S}$ of the MLSVD distributes the ``energy'' (i.e., the magnitude of its entries) in such a way so that it concentrates more energy at the first entry $s_{11 \ldots 1}$ and disperses as we move along each dimension. Figure~\ref{fig:slices-energy} illustrates the energy distribution when $\mathcal{S}$ is a third order tensor. The red slices contains more energy and it changes to white when the slice contains less energy. Note that the energy of the slices are given precisely by the singular values.
				
	\begin{definition}
		The \emph{multilinear rank of $\mathcal{T}$} is the $L$-tuple $\left( R_1, \ldots, R_L \right)$. We denote $rank_\boxplus(\mathcal{T}) = (R_1, \ldots, R_L)$.
	\end{definition}
	
	The multilinear rank is also often called the \emph{Tucker rank} of $\mathcal{T}$. Now we will see some important results regarding this rank. More details can be found in \cite{lim}.
			
	\begin{theorem}[V. de Silva, and L.H. Lim, 2008]
		Let $\mathcal{T} \in \mathbb{R}^{I_1} \otimes \ldots \otimes \mathbb{R}^{I_L}$ and $\textbf{M}^{(1)} \in \mathbb{R}^{I_1' \times I_1}, \ldots, \textbf{M}^{(L)} \in \mathbb{R}^{I_L' \times I_L}$. Then the following statements holds.
		\begin{enumerate}
			\item $\| rank_\boxplus(\mathcal{T}) \|_\infty \leq rank(T)$.
			\item $rank_\boxplus((\textbf{M}^{(1)}, \ldots, \textbf{M}^{(L)})\cdot \mathcal{T}) \leq rank_\boxplus(\mathcal{T})$.
			\item If $\textbf{M}^{(1)} \in GL(I_1'), \ldots, \textbf{M}^{(L)} \in GL(I_L')$, then $rank_\boxplus((\textbf{M}^{(1)}, \ldots, \textbf{M}^{(L)}) \cdot \mathcal{T}) = rank_\boxplus(\mathcal{T})$, where $GL(d)$ denotes the group of invertible matrices $d \times d$.
		\end{enumerate}
	\end{theorem}
	
	After computing the MLSVD of $\mathcal{T}$, notice that computing a CPD for $\mathcal{S}$ is equivalent to computing a CPD for $\mathcal{T}$. Indeed, if $\mathcal{S} = \displaystyle \sum_{r=1}^R \textbf{w}_r^{(1)} \otimes \ldots \otimes \textbf{w}_r^{(L)}$, then we can write $\mathcal{S} = (\textbf{W}^{(1)}, \ldots, \textbf{W}^{(L)}) \cdot \mathcal{I}_{R \times \ldots \times R}$, which implies that
	$$\mathcal{T} = (\textbf{U}^{(1)}, \ldots, \textbf{U}^{(L)}) \cdot \left( (\textbf{W}^{(1)}, \ldots, \textbf{W}^{(L)}) \cdot \mathcal{I}_{R \times \ldots \times R} \right) = $$
	$$ = (\textbf{U}^{(1)} \textbf{W}^{(1)}, \ldots, \textbf{U}^{(L)} \textbf{W}^{(L)}) \cdot \mathcal{I}_{R \times \ldots \times R}.$$
			
\begin{figure}[htbp]
		\centering
		\includegraphics[scale=.45]{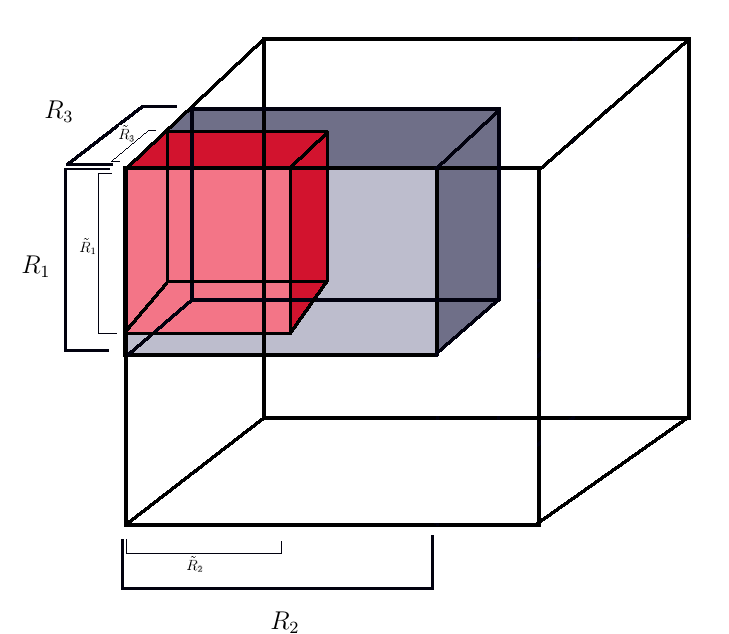}
		\caption{Truncated tensor $\tilde{\mathcal{S}}$.}
		\label{fig:trunc}
	\end{figure}

	Let $rank_\boxplus(\mathcal{T}) = (R_1, \ldots, R_L)$. We use the algorithm introduced by N. Vannieuwenhove, R. Vandebril and K. Meerberge \cite{sthosvd} to compute a truncated MLSVD. Now let $(\tilde{R}_1, \ldots, \tilde{R}_L) \leq (R_1, \ldots, R_L)$ be a lower multilinear rank. We define $\tilde{\textbf{U}}^{(\ell)} = \left[ \textbf{U}_{:1}^{(\ell)}, \ldots, \textbf{U}_{:\tilde{R}_\ell}^{(\ell)} \right] \in \mathbb{R}^{I_\ell \times \tilde{R}_\ell}$ to be the matrix composed by the first columns of $\textbf{U}^{(\ell)}$, and $\tilde{\mathcal{S}} \in \mathbb{R}^{\tilde{R}_1 \times \ldots \times \tilde{R}_L}$ is such that $\tilde{s}_{i_1 \ldots i_L} = s_{i_1 \ldots i_L}$ for $1 \leq i_1 \leq \tilde{R}_1, \ldots, 1 \leq i_L \leq \tilde{R}_L$. Figure~\ref{fig:trunc} illustrates such a truncation in the case of a third order tensor. The white part correspond to $\mathcal{S}$ after we computed the full MLSVD, the gray tensor is the reduced format of $\mathcal{S}$, and the red tensor is the truncated tensor $\tilde{\mathcal{S}}$. Our goal is to find the smallest $(\tilde{R}_1, \ldots, \tilde{R}_L)$ such that $\| \mathcal{S} - \tilde{\mathcal{S}} \|$ is not so large\footnote{Actually, $\mathcal{S}$ and $\tilde{\mathcal{S}}$ belongs to different spaces. We committed an abuse of notation and wrote $\| \mathcal{S} - \tilde{\mathcal{S}} \|$ considering the projection of $\tilde{\mathcal{S}}$ over the space of $\mathcal{S}$, that is, enlarge $\tilde{\mathcal{S}}$ so it has the same size of $\mathcal{S}$ and consider these new entries as zeros.}. Since reducing $(\tilde{R}_1, \ldots, \tilde{R}_L)$ too much causes $\| \mathcal{S} - \tilde{\mathcal{S}} \|$ to increase, there is a trade off we have to manage in the best way possible. 
			
\section{Algorithm description and analysis}
\label{sec:alg}

Given a tensor $\mathcal{T} \in \mathbb{R}^{m \times n \times p}$ and a rank value $R$, we want to find a rank $R$ approximation for $\mathcal{T}$. More precisely, we want to find vectors $\textbf{x}_1, \ldots, \textbf{x}_R \in \mathbb{R}^m$, $\textbf{y}_1, \ldots, \textbf{y}_R \in \mathbb{R}^n$ and $\textbf{z}_1, \ldots, \textbf{z}_R \in \mathbb{R}^p$ such that
	\begin{equation}
		\mathcal{T} \approx \sum_{r=1}^R \textbf{x}_r \otimes \textbf{y}_r \otimes \textbf{z}_r. \label{eq:rank_approx}
	\end{equation}			
	
	Finding the best approximation is known to be a ill-posed problem \cite{lim} so we will be content with just a reasonable approximation. This also depends on the choice of $R$ but we won't address this issue here. For each problem at hand one may have a good guess of what should be the rank. In the worse scenario it is still possible to try several ranks and keep the best fit, although caution is necessary to not overfit the data.
	
	Consider the factors $\textbf{X} = [ \textbf{x}_1, \ldots, \textbf{x}_R ] \in \mathbb{R}^{m \times R}$, $\textbf{Y} = [ \textbf{y}_1, \ldots, \textbf{y}_R ] \in \mathbb{R}^{n \times R}$, $\textbf{Z} = [ \textbf{z}_1, \ldots, \textbf{z}_R ] \in \mathbb{R}^{p \times R}$. Using the multilinear multiplication notation we can write ~\ref{eq:rank_approx} as $\mathcal{T} \approx (\textbf{X}, \textbf{Y}, \textbf{Z}) \cdot \mathcal{I}_{R \times R \times R}$. This view may lead us to consider minimizing the map 
	\begin{equation}
	(\textbf{X}, \textbf{Y}, \textbf{Z}) \mapsto \frac{1}{2} \left\| \mathcal{T} - (\textbf{X}, \textbf{Y}, \textbf{Z}) \cdot \mathcal{I}_{R \times R \times R} \right\|^2. \label{eq:rank_approx_tucker}
	\end{equation}
	
	Since the matrix structure is not really useful in~\ref{eq:rank_approx_tucker}, we may just vectorize the matrices and concatenate them all to form the single vector
	$$\textbf{w} = 
	\left[ 
	\begin{array}{c}
		vec(\textbf{X})\\
		vec(\textbf{Y})\\
		vec(\textbf{Z})
	\end{array}
	\right].$$
	The map $vec$ just stacks all columns of a matrix, with the first column at the top and going on until the last column at the bottom. With this notation we consider the error function $F:\mathbb{R}^{R(m+n+p)} \to \mathbb{R}$ given by
	$$F(\textbf{w}) = \frac{1}{2} \left\| \mathcal{T} - (\textbf{X}, \textbf{Y}, \textbf{Z}) \cdot \mathcal{I}_{R \times R \times R} \right\|^2 = $$
	
	$$ = \frac{1}{2} \left\| \mathcal{T} - \sum_{r=1}^R \textbf{x}_r \otimes \textbf{y}_r \otimes \textbf{z}_r \right\|^2 = $$	
	
	$$ = \frac{1}{2} \sum_{i=1}^m \sum_{j=1}^n \sum_{k=1}^p \left( \mathcal{T}_{ijk} - \sum_{r=1}^R \textbf{x}_{ir} \textbf{y}_{jr} \textbf{z}_{kr} \right)^2.$$	
	
	For each choice of $\textbf{w}$ we write the approximating tensor as $\tilde{\mathcal{T}} = \tilde{\mathcal{T}}(\textbf{w})$. In order to minimize $F$ we first introduce the Gauss-Newton algorithm (without damping). Each difference $\mathcal{T}_{ijk} - \tilde{\mathcal{T}}_{ijk} = \mathcal{T}_{ijk} - \sum_{r=1}^R \textbf{x}_i^{(r)} \textbf{y}_j^{(r)} \textbf{z}_k^{(r)}$ is the residual between $\mathcal{T}_{ijk}$ and $\tilde{\mathcal{T}}_{ijk}$.	From these residuals we form the function $f:\mathbb{R}^{R(m+n+p)} \to \mathbb{R}^{mnp}$ defined by $f = (f_{111}, f_{112}, \ldots, f_{mnp})$, where $f_{ijk}(\textbf{w}) = \mathcal{T}_{ijk} - \tilde{\mathcal{T}}_{ijk}$. Then we can write
	$$F(\textbf{w}) = \frac{1}{2} \sum_{i=1}^m \sum_{j=1}^n \sum_{k=1}^p f_{ijk}(\textbf{w})^2 = \frac{1}{2} \| f(\textbf{w}) \|^2.$$	
	
	\begin{lemma} The partial derivatives of $f$ are the following.
		$$\frac{\partial f_{ijk}}{\partial x_{i'r}} = \left\{
		\begin{array}{cc}
    		- y_{jr} z_{kr}, & \text{if } i' = i\\
    		0, & \text{otherwise}
		\end{array}
		\right.$$

		$$\frac{\partial f_{ijk}}{\partial y_{j'r}} = \left\{
		\begin{array}{cc}
    		- x_{ir} z_{kr}, & \text{if } j' = j\\
    		0, & \text{otherwise}
		\end{array}
		\right.$$

		$$\frac{\partial f_{ijk}}{\partial z_{k'r}} = \left\{
		\begin{array}{cc}
    		- x_{ir} y_{jr}, & \text{if } k' = k\\
    		0, & \text{otherwise}
		\end{array}
		\right.$$
	\end{lemma}
	
	Denote the derivative of $f$ at $\textbf{w}$ (the Jacobian matrix) by $\textbf{J}_f = \textbf{J}_f(\textbf{w})$. Also denote the gradient of $F$ at $\textbf{w}$ by $\nabla F = \nabla F(\textbf{w})$ and the corresponding Hessian by $\textbf{H}_F = \textbf{H}_F(\textbf{w})$. We can relate the derivatives of $f$ and $F$ through the following result.
	
	\begin{lemma} \label{relations} The following identities always holds.
		\begin{enumerate}
			\item $\displaystyle \nabla F = \textbf{J}_f^T \cdot f$.
			\item $\displaystyle \textbf{H}_F = \textbf{J}_f^T \textbf{J}_f + \sum_{i=1}^m \sum_{j=1}^n \sum_{k=1}^p f_{ijk} \cdot \textbf{H}_{f_{ijk}}$, where $\textbf{H}_{f_{ijk}}$ is the Hessian matrix of $f_{ijk}$.
		\end{enumerate}
	\end{lemma}
	
	As the algorithm converges we expect to have $f_{ijk} \approx 0$ for all $i,j,k$. Together with lemma ~\ref{relations} this observation shows that $\textbf{H}_F \approx \textbf{J}_f^T \textbf{J}_f$ when close to a optimal point.
	
	\subsection{Gauss-Newton algorithm}	
		The Gauss-Newton algorithm is obtained by first considering a first order approximation of $f$ at a point $\textbf{w}^{(0)} \in \mathbb{R}^{R(m+n+p)}$, that is,
		\begin{equation}
		f(\textbf{w}^{(0)} + \underbrace{(\textbf{w} - \textbf{w}^{(0)})}_{step} ) = f(\textbf{w}) \approx f(\textbf{w}^{(0)}) + \textbf{J}_f(\textbf{w}^{(0)}) \cdot (\textbf{w} - \textbf{w}^{(0)}). \label{eq:linear_approx}
		\end{equation}
	
		In order to minimize ~\ref{eq:linear_approx} at the neighborhood of $\textbf{w}^{(0)}$ we can compute the minimum of $\| f(\textbf{w}^{(0)}) + \textbf{J}_f(\textbf{w}^{(0)}) \cdot (\textbf{w} - \textbf{w}^{(0)}) \|$ for $\textbf{w} \in \mathbb{R}^{R(m+n+p)}$.	 Note that minimizing ~\ref{eq:linear_approx} is a least squares problem since we can rewrite this problem as $\displaystyle \min_\textbf{x} \| \textbf{A}\textbf{x} - \textbf{b}\|$ for $\textbf{A} = \textbf{J}_f(\textbf{w}^{(0)})$, $\textbf{x} = \textbf{w} - \textbf{w}^{(0)}$, $\textbf{b} = - f(\textbf{w}^{(0)})$.
	
		The solution gives us the next iterate $\textbf{w}^{(1)}$. More generally, we obtain $\textbf{w}^{(k+1)}$ from $\textbf{w}^{(k)}$ by defining $\textbf{w}^{(k+1)} = \textbf{x}^\ast + \textbf{w}^{(k)}$, where $\textbf{x}^\ast$ is the solution of the normal equations
		\begin{equation}
			\textbf{A}^T \textbf{A} \textbf{x} = \textbf{A}^T \textbf{b} \label{eq:normal_eq}
		\end{equation}
for $\textbf{A} = \textbf{J}_f(\textbf{w}^{(k)})$, $\textbf{x} = \textbf{w} - \textbf{w}^{(k)}$, $\textbf{b} = - f(\textbf{w}^{(k)})$. This iteration process is the Gauss-Newton algorithm, and it is guaranteed to converge to a local minimum \cite{madsen}.

		The Jacobian matrix $\textbf{J}_f$ is always singular so the problem is ill-posed. For this reason we want to regularize the problem. By introducing regularization we can avoid singularity and improve convergence. A common approach is to introduce a suitable regularization matrix $\textbf{L} \in \mathbb{R}^{R(m+n+p) \times R(m+n+p)}$ called \emph{Tikhonov matrix}. Instead of solving equation ~\ref{eq:normal_eq} we solve 
		\begin{equation}
			\left( \textbf{A}^T \textbf{A} + \textbf{L}^T \textbf{L} \right) \textbf{x} = \textbf{A}^T \textbf{b}. \label{eq:reg_eq}
		\end{equation}
	
		When $\textbf{L}$ is diagonal, this is called the dGN algorithm. TFX works with a Tikhonov matrix of the form $\textbf{L} = \mu \textbf{D}$, where $\textbf{D}$ is a certain $R(m+n+p) \times R(m+n+p)$ matrix with positive diagonal and $\mu > 0$ is the \emph{damping parameter}. The important property of $\textbf{D}$ is that $\textbf{A}^T \textbf{A} + \textbf{D}$ is diagonally dominant. The damping parameter $\mu$ is usually updated at each iteration. These updates are very important since $\mu$ influences both the direction and the size of the step at each iteration. Let $\tilde{\mathcal{T}}^{(k)}$ be the approximating tensor at the $k$-th iteration and $\tilde{f}(\textbf{w}^{(k)}) = f(\textbf{w}^{(k-1)}) + \textbf{J}_f(\textbf{w}^{(k-1)}) \cdot (\textbf{w}^{(k)} - \textbf{w}^{(k-1)})$ is the first order approximation of $f$ at $\textbf{w}^{(k)}$. TFX uses the update strategy
		
		\begin{flalign*}
			& \texttt{if } g < 0.75\\
			& \hspace{.8cm} \mu \leftarrow \mu/2\\
			& \texttt{else if } g > 0.9\\
			& \hspace{.8cm} \mu \leftarrow 1.5 \cdot \mu
		\end{flalign*}
where $g$ is the \emph{gain ratio}, defined as
		$$g = \frac{ \| \mathcal{T} - \tilde{\mathcal{T}}^{(k-1)} \|^2 - \| \mathcal{T} - \tilde{\mathcal{T}}^{(k)} \|^2 }{ \| \mathcal{T} - \tilde{\mathcal{T}}^{(k-1)} \|^2 - \| \tilde{f}(\textbf{w}^{(k)}) \|^2 }.$$
The denominator is the predicted improvement from the $(k-1)$-th iteration to $k$-th iteration, whereas the numerator measures the actual improvement. 
	
		\begin{theorem} With the notations above, the following holds.
			\begin{enumerate}
				\item $\textbf{J}_f(\textbf{w}^{(k)}) \textbf{J}_f(\textbf{w}^{(k)}) + \mu \textbf{D}$ is a positive definite matrix for all $\mu > 0$.
				\item $\textbf{w}^{(k+1)} - \textbf{w}^{(k)}$ is a descent direction for $F$ at $\textbf{w}^{(k)}$.
				\item If $\mu$ is big enough, then $\small \textbf{w}^{(k+1)} - \textbf{w}^{(k)} \approx - \frac{1}{\mu} \textbf{J}_f(\textbf{w}^{(k)})^T \cdot f(\textbf{w}^{(k)}) = - \frac{1}{\mu} \nabla F(\textbf{w}^{(k)})$.
				\item If $\mu$ is small enough, then $\textbf{w}^{(k+1)} - \textbf{w}^{(k)} \approx \textbf{w}_{GN}^{(k+1)} - \textbf{w}^{(k)}$, where $\textbf{w}_{GN}^{(k+1)}$ is the point we would obtain using classic Gauss-Newton iteration (i.e.,  without regularization). 
			\end{enumerate}
		\end{theorem} 
	
		\begin{remark} 
			Item 3 is to be used when the current iteration is far from the solution, since $ - \frac{1}{\mu} \nabla F(\textbf{w}^{(k)})$ is a short step in the descent direction. We want to be careful when distant to the solution. This shows that dGN behaves as the gradient descent algorithm when distant to the solution. On the other hand, item 4 is to be used when the current iteration is close to the solution, since the step is closer to the classic Gauss-Newton, we may attain quadratic convergence at the final.
		\end{remark}
		
	\subsection{Exploiting the structure of $\textbf{A}^T \textbf{A}$}	
		At each iteration we must solve equation ~\ref{eq:reg_eq}, and since $\textbf{A}^T \textbf{A} + \textbf{L}^T \textbf{L}$ is positive definite we are able to find an approximated solution by taking just a few iterations of the conjugate gradient algorithm. However $\textbf{A}^T \textbf{A} + \textbf{L}^T \textbf{L}$ is a $R(m+n+p) \times R(m+n+p)$ matrix, and at each iteration this may be costly since each matrix-vector  multiplication costs $\mathcal{O}(R^2(m+n+p)^2)$ flops (floating point operations). We can work around this issue by exploiting the structure of $\textbf{A}^T \textbf{A}$. 
		
		\begin{theorem}
			Let $\textbf{A} = \textbf{J}_f(\textbf{w}^{(k)})$ and let $\textbf{X}, \textbf{Y}, \textbf{Z}$ be the factor matrices of the approximating CPD at $k$-th iteration, where $\textbf{X} = [ \textbf{x}^{(1)}, \ldots, \textbf{x}^{(R)} ]$, $\textbf{Y} = [ \textbf{y}^{(1)}, \ldots, \textbf{y}^{(R)} ]$, $\textbf{Z} = [ \textbf{z}^{(1)}, \ldots, \textbf{z}^{(R)} ]$. Then,		
	
		$$\textbf{A}^T \textbf{A} = 
	 	\left[ \begin{array}{ccc}
	 		\textbf{B}_X & \textbf{B}_{XY} & \textbf{B}_{XZ}\\
	 		\textbf{B}_{XY}^T & \textbf{B}_Y & \textbf{B}_{YZ}\\
	 		\textbf{B}_{XZ}^T & \textbf{B}_{YZ}^T & \textbf{B}_Z
	 	\end{array} \right],$$
where
	
		$$\textbf{B}_X = 
	 	\left[ \begin{array}{ccc}
	 		\langle \textbf{y}_1, \textbf{y}_1 \rangle \langle \textbf{z}_1, \textbf{z}_1 \rangle \textbf{I}_m & \ldots & \langle \textbf{y}_1, \textbf{y}_R \rangle \langle \textbf{z}_1, \textbf{z}_R \rangle \textbf{I}_m\\
	 		\vdots & & \vdots\\
	 		\langle \textbf{y}_R, \textbf{y}_1 \rangle \langle \textbf{z}_R, \textbf{z}_1 \rangle \textbf{I}_m & \ldots & \langle \textbf{y}_R, \textbf{y}_R \rangle \langle \textbf{z}_R, \textbf{z}_R \rangle \textbf{I}_m\\
	 	\end{array} \right],$$
	 
	 	$$\textbf{B}_Y = 
	 	\left[ \begin{array}{ccc}
	 		\langle \textbf{x}_1, \textbf{x}_1 \rangle \langle \textbf{z}_1, \textbf{z}_1 \rangle \textbf{I}_n & \ldots & \langle \textbf{x}_1, \textbf{x}_R \rangle \langle \textbf{z}_1, \textbf{z}_R \rangle \textbf{I}_n\\
	 		\vdots & & \vdots\\
	 		\langle \textbf{x}_R, \textbf{x}_1 \rangle \langle \textbf{z}_R, \textbf{z}_1 \rangle \textbf{I}_n & \ldots & \langle \textbf{x}_R, \textbf{x}_R \rangle \langle \textbf{z}_R, \textbf{z}_R \rangle \textbf{I}_n\\
	 	\end{array} \right],$$
	 
	 	$$\textbf{B}_Z = 
	 	\left[ \begin{array}{ccc}
	 		\langle \textbf{x}_1, \textbf{x}_1 \rangle \langle \textbf{y}_1, \textbf{y}_1 \rangle \textbf{I}_p & \ldots & \langle \textbf{y}_1, \textbf{y}_R \rangle \langle \textbf{y}_1, \textbf{y}^{(r)} \rangle \textbf{I}_p\\
	 		\vdots & & \vdots\\
	 		\langle \textbf{x}_R, \textbf{x}_1 \rangle \langle \textbf{y}_R, \textbf{y}_1 \rangle \textbf{I}_p & \ldots & \langle \textbf{x}_R, \textbf{x}_R \rangle \langle \textbf{y}_R, \textbf{y}_R \rangle \textbf{I}_p\\
	 	\end{array} \right],$$
	 
	 	$$\textbf{B}_{XY} = 
	 	\left[ \begin{array}{ccc}
	 		\langle \textbf{z}_1, \textbf{z}_1 \rangle \textbf{x}_1 \textbf{y}_1^T & \ldots & \langle \textbf{z}_1, \textbf{z}_R \rangle \textbf{x}^{(r)} \textbf{y}_1^T\\
	 		\vdots & & \vdots\\
	 		\langle \textbf{z}_R, \textbf{z}_1 \rangle \textbf{z}_1 \textbf{y}_R^T & \ldots & \langle \textbf{z}_R, \textbf{z}_R \rangle \textbf{x}_R \textbf{y}_R^T\\
	 	\end{array} \right],$$
	 
	 	$$\textbf{B}_{XZ} = 
	 	\left[ \begin{array}{ccc}
	 		\langle \textbf{y}_1, \textbf{y}_1 \rangle \textbf{x}_1 \textbf{z}_1^T & \ldots & \langle \textbf{y}_1, \textbf{y}_R \rangle \textbf{x}_R (\textbf{z}^{(1)})^T\\
	 		\vdots & & \vdots\\
	 		\langle \textbf{y}_R, \textbf{y}_1 \rangle \textbf{x}_1 \textbf{z}_R^T & \ldots & \langle \textbf{y}_R, \textbf{y}_R \rangle \textbf{x}_R \textbf{z}_R^T\\
	 	\end{array} \right],$$
	 
	 	$$\textbf{B}_{YZ} = 
	 	\left[ \begin{array}{ccc}
	 		\langle \textbf{x}_1, \textbf{x}_1 \rangle \textbf{y}_1 \textbf{z}_1^T & \ldots & \langle \textbf{x}_1, \textbf{x}^{(r)} \rangle \textbf{y}^{(r)} \textbf{z}_1^T\\
	 		\vdots & & \vdots\\
	 		\langle \textbf{x}_R, \textbf{x}_1 \rangle \textbf{y}_1 \textbf{z}_R^T & \ldots & \langle \textbf{x}_R, \textbf{x}_R \rangle \textbf{y}_R \textbf{z}_R^T\\
	 	\end{array} \right].$$
	 	
	 \end{theorem}
	 
	 Now we can see how to retrieve $\textbf{A}^T \textbf{A}$ from the factor matrices $\textbf{X}, \textbf{Y}, \textbf{Z}$ with few computations and low memory cost. First we compute and store all the scalar products. It is easy to see we will need to store $9 R^2$ floats. The other parts of $\textbf{A}^T \textbf{A}$ can be obtained directly from the factor matrices so we are done with regard to memory costs. This is a big reduction in memory size since the original matrix would require $R^2 \left( m+n+p \right)^2$ floats in dense format. The overall cost to compute those scalar products is $\mathcal{O}\left( R^2 (m+n+p) \right)$ flops, which is also reasonable.
	 
	 It is convenient to store the above products in matrix form, so we write
			
		$$\Pi_X = 
	 	\left[ \begin{array}{ccc}
	 		\langle \textbf{y}_1, \textbf{y}_1 \rangle \langle \textbf{z}_1, \textbf{z}_1 \rangle & \ldots & \langle \textbf{y}_1, \textbf{y}_R \rangle \langle \textbf{z}_1, \textbf{z}_R \rangle \\
	 		\vdots & & \vdots\\
	 		\langle \textbf{y}_R, \textbf{y}_1 \rangle \langle \textbf{z}_R, \textbf{z}_1 \rangle & \ldots & \langle \textbf{y}_R, \textbf{y}_R \rangle \langle \textbf{z}_R, \textbf{z}_R \rangle \\
	 	\end{array} \right],$$
	 
	 	$$\Pi_Y = 
	 	\left[ \begin{array}{ccc}
	 		\langle \textbf{x}_1, \textbf{x}_1 \rangle \langle \textbf{z}_1, \textbf{z}_1 \rangle & \ldots & \langle \textbf{x}_1, \textbf{x}_R \rangle \langle \textbf{z}_1, \textbf{z}_R \rangle \\
	 		\vdots & & \vdots\\
	 		\langle \textbf{x}_R, \textbf{x}_1 \rangle \langle \textbf{z}_R, \textbf{z}_1 \rangle & \ldots & \langle \textbf{x}_R, \textbf{x}_R \rangle \langle \textbf{z}_R, \textbf{z}_R \rangle \\
	 	\end{array} \right],$$
	 
	 	$$\Pi_Z = 
	 	\left[ \begin{array}{ccc}
	 		\langle \textbf{x}_1, \textbf{x}_1 \rangle \langle \textbf{y}_1, \textbf{y}_1 \rangle & \ldots & \langle \textbf{y}_1, \textbf{y}_R \rangle \langle \textbf{y}_1, \textbf{y}^{(r)} \rangle\\
	 		\vdots & & \vdots\\
	 		\langle \textbf{x}_R, \textbf{x}_1 \rangle \langle \textbf{y}_R, \textbf{y}_1 \rangle & \ldots & \langle \textbf{x}_R, \textbf{x}_R \rangle \langle \textbf{y}_R, \textbf{y}_R \rangle \\
	 	\end{array} \right],$$
	 
	 	$$\Pi_{XY} = 
	 	\left[ \begin{array}{ccc}
	 		\langle \textbf{z}_1, \textbf{z}_1 \rangle & \ldots & \langle \textbf{z}_1, \textbf{z}_R \rangle \\
	 		\vdots & & \vdots\\
	 		\langle \textbf{z}_R, \textbf{z}_1 \rangle & \ldots & \langle \textbf{z}_R, \textbf{z}_R \rangle \\
	 	\end{array} \right],$$
	 
	 	$$\Pi_{XZ} = 
	 	\left[ \begin{array}{ccc}
	 		\langle \textbf{y}_1, \textbf{y}_1 \rangle & \ldots & \langle \textbf{y}_1, \textbf{y}_R \rangle \\
	 		\vdots & & \vdots\\
	 		\langle \textbf{y}_R, \textbf{y}_1 \rangle & \ldots & \langle \textbf{y}_R, \textbf{y}_R \rangle \\
	 	\end{array} \right],$$
	 
	 	$$\Pi_{YZ} = 
	 	\left[ \begin{array}{ccc}
	 		\langle \textbf{x}_1, \textbf{x}_1 \rangle & \ldots & \langle \textbf{x}_1, \textbf{x}^{(r)} \rangle \\
	 		\vdots & & \vdots\\
	 		\langle \textbf{x}_R, \textbf{x}_1 \rangle & \ldots & \langle \textbf{x}_R, \textbf{x}_R \rangle \\
	 	\end{array} \right].$$

		As already mentioned, $\textbf{A}^T \textbf{A}$ will be used to solve the normal equations~\ref{eq:reg_eq}. The algorithm of choice to accomplish this is the conjugate gradient method. This classical algorithm is particularly efficient to solve normal equations where the matrix is positive definite, which is our case. Furthermore, our version of the conjugate gradient is \emph{matrix-free}, that is, we are able to compute matrix-vector products $\textbf{A}^T \textbf{A} \cdot \textbf{v}$ without actually constructing $\textbf{A}^T \textbf{A}$. By exploiting the block structure of $\textbf{A}^T \textbf{A}$ we can save memory and the computational cost still is lower than the naive cost of $R^2 \left( m+n+p \right)^2$ flops.   
		
		\begin{theorem} \label{Hv}
			Given any vector $\textbf{v} \in \mathbb{R}^{R (m+n+p)}$, write 			
			$$\textbf{v} = 
			\left[ \begin{array}{c}
				vec(\textbf{V}_X)\\
				vec(\textbf{V}_Y)\\
				vec(\textbf{V}_Z)
			\end{array} \right]$$
where $\textbf{V}_X \in \mathbb{R}^{m \times R}$, $\textbf{V}_Y \in \mathbb{R}^{n \times R}$, $\textbf{V}_Z \in \mathbb{R}^{p \times R}$. Then, 

			$$\textbf{A}^T \textbf{A} \cdot \textbf{v} = \footnotesize \left[ \begin{array}{c}
					\displaystyle vec\left( \textbf{V}_X \cdot \Pi_X \right) + vec\left( \textbf{X} \cdot \Big( \Pi_{XY} \ast \big( \textbf{V}_Y^T \cdot \textbf{Y} \big) \Big) \right) + vec\left( \textbf{X} \cdot \Big( \Pi_{XZ} \ast \big( \textbf{V}_Z^T \cdot \textbf{Z} \big) \Big) \right)\\ \\
					\displaystyle vec\left( \textbf{Y} \cdot \Big( \Pi_{XY} \ast \big( \textbf{V}_X^T \cdot \textbf{X} \big) \Big) \right) + vec\left( \textbf{V}_Y \cdot \Pi_Y \right) + vec\left( \textbf{Y} \cdot \Big( \Pi_{YZ} \ast \big( \textbf{V}_Z^T \cdot \textbf{Z} \big) \Big) \right)\\ \\
					\displaystyle vec\left( \textbf{Z} \cdot \Big( \Pi_{XZ} \ast \big( \textbf{V}_X^T \cdot \textbf{X} \big) \Big) \right) + vec\left( \textbf{Z} \cdot \Big( \Pi_{YZ} \ast \big( \textbf{V}_Y^T \cdot \textbf{Y} \big) \Big) \right) + vec\left( \textbf{V}_Z \cdot \Pi_Z \right)
				\end{array} \right].$$		
		\end{theorem}
		\bigskip
		
		\textbf{Proof:} We will prove only the equality for the first block 
		$$\displaystyle vec\left( \textbf{V}_X \ast \Pi_X \right) + vec\left( \textbf{X} \cdot \Big( \Pi_{XY} \ast \big( \textbf{V}_Y^T \cdot \textbf{Y} \big) \Big) \right) + vec\left( \textbf{X} \cdot \Big( \Pi_{XZ} \ast \big( \textbf{V}_Z^T \cdot \textbf{Z} \big) \Big) \right),$$ 
the others are analogous. First notice that

			$$\textbf{A}^T \textbf{A} \cdot \textbf{v} = \textbf{J}_f^T \textbf{J}_f \cdot \textbf{v} = 
			\left[ \begin{array}{ccc}
			 	\displaystyle\frac{\partial f}{\partial \textbf{X}}^T \frac{\partial f}{\partial \textbf{X}} \ & \ \displaystyle\frac{\partial f}{\partial \textbf{X}}^T \frac{\partial f}{\partial \textbf{Y}} \ & \ \displaystyle\frac{\partial f}{\partial \textbf{X}}^T \frac{\partial f}{\partial \textbf{Z}}\\ \\
			  	\displaystyle\frac{\partial f}{\partial \textbf{Y}}^T \frac{\partial f}{\partial \textbf{X}} \ & \ \displaystyle\frac{\partial f}{\partial \textbf{Y}}^T \frac{\partial f}{\partial \textbf{Y}} \ & \ \displaystyle\frac{\partial f}{\partial \textbf{Y}}^T \frac{\partial f}{\partial \textbf{Z}}\\ \\
			   \displaystyle\frac{\partial f}{\partial \textbf{Z}}^T \frac{\partial f}{\partial \textbf{X}} \ & \ \displaystyle\frac{\partial f}{\partial \textbf{Z}}^T \frac{\partial f}{\partial \textbf{Y}} \ & \ \displaystyle\frac{\partial f}{\partial \textbf{Z}}^T \frac{\partial f}{\partial \textbf{Z}}
			\end{array} \right]
			\left[ \begin{array}{c}
				vec(\textbf{V}_X)\\
				vec(\textbf{V}_Y)\\
				vec(\textbf{V}_Z)
			\end{array} \right] = $$
			
			$$\left[ \begin{array}{c}
			 	\displaystyle\frac{\partial f}{\partial \textbf{X}}^T \frac{\partial f}{\partial \textbf{X}} \cdot vec(\textbf{V}_X) + \displaystyle\frac{\partial f}{\partial \textbf{X}}^T \frac{\partial f}{\partial \textbf{Y}} \cdot vec(\textbf{V}_Y) + \displaystyle\frac{\partial f}{\partial \textbf{X}}^T \frac{\partial f}{\partial \textbf{Z}} \cdot vec(\textbf{V}_Z)\\ \\
			  	\displaystyle\frac{\partial f}{\partial \textbf{Y}}^T \frac{\partial f}{\partial \textbf{X}} \cdot vec(\textbf{V}_X) + \displaystyle\frac{\partial f}{\partial \textbf{Y}}^T \frac{\partial f}{\partial \textbf{Y}} \cdot vec(\textbf{V}_Y) + \displaystyle\frac{\partial f}{\partial \textbf{Y}}^T \frac{\partial f}{\partial \textbf{Z}} \cdot vec(\textbf{V}_Z)\\ \\
			   \displaystyle\frac{\partial f}{\partial \textbf{Z}}^T \frac{\partial f}{\partial \textbf{X}} \cdot vec(\textbf{V}_X) + \displaystyle\frac{\partial f}{\partial \textbf{Z}}^T \frac{\partial f}{\partial \textbf{Y}} \cdot vec(\textbf{V}_Y) + \displaystyle\frac{\partial f}{\partial \textbf{Z}}^T \frac{\partial f}{\partial \textbf{Z}} \cdot vec(\textbf{V}_Z)
			\end{array} \right], $$
			
where

		$$\frac{\partial f}{\partial \textbf{X}} = \left[
		\begin{array}{ccccccc}
			\displaystyle\frac{\partial f_{111}}{partial \textbf{X}_{11}} & \ldots & \displaystyle\frac{\partial f_{111}}{\partial \textbf{X}_{1R}} & \ldots & \displaystyle\frac{\partial f_{111}}{\partial \textbf{X}_{m1}} & \ldots & \displaystyle\frac{\partial f_{111}}{\partial \textbf{X}_{mR}} \\ \\
			\displaystyle\frac{\partial f_{112}}{\partial \textbf{X}_{11}} & \ldots & \displaystyle\frac{\partial f_{112}}{\partial \textbf{X}_{1R}} & \ldots & \displaystyle\frac{\partial f_{112}}{\partial \textbf{X}_{m1}} & \ldots & \displaystyle\frac{\partial f_{112}}{\partial \textbf{X}_{mR}} \\ \\
			\vdots & & \vdots & & \vdots & & \vdots\\ \\
			\displaystyle\frac{\partial f_{mnp}}{\partial \textbf{X}_{11}} & \ldots & \displaystyle\frac{\partial f_{mnp}}{\partial \textbf{X}_{1R}} & \ldots & \displaystyle\frac{\partial f_{mnp}}{\partial \textbf{X}_{m1}} & \ldots & \displaystyle\frac{\partial f_{mnp}}{\partial \textbf{X}_{mR}}
		\end{array}
		\right],$$
and the other derivatives are defined similarly.

		Now we simplify each term in the summation above. It is necessary to consider two separate cases. Again, we only prove one particular case since the other ones are proven similarly.\bigskip
			
			\textbf{Case 1 (different modes):} Write $\textbf{V}_{Y} = [\textbf{v}_{Y_1}, \ldots, \textbf{v}_{Y_R}]$, where each $\textbf{v}_{Y_r} \in \mathbb{R}^n$ is a column of $\textbf{V}_{Y}$. Then
			
			$$\frac{\partial f}{\partial \textbf{X}}^T \frac{\partial f}{\partial \textbf{Y}} \cdot vec(\textbf{V}_Y) = $$
			
			$$ = 
			\small\left[ \begin{array}{ccc}
				\displaystyle \langle \textbf{z}_1, \textbf{z}_1 \rangle \cdot \textbf{x}_1 \textbf{y}_1^T & \ldots & \displaystyle \langle \textbf{z}_1, \textbf{z}_R \rangle \cdot \textbf{x}_R \textbf{y}_1^T\\
				\vdots & & \vdots\\
				\displaystyle \langle \textbf{z}_R, \textbf{z}_1 \rangle \cdot \textbf{x}_1 \textbf{y}_R^T & \ldots & \displaystyle \langle \textbf{z}_R, \textbf{z}_R \rangle \cdot \textbf{x}_R \textbf{y}_R^T\\
			\end{array} \right] 
			\left[ \begin{array}{c}
				\textbf{v}_{Y_1}\\
				\vdots\\
				\textbf{v}_{Y_R}
			\end{array} \right] = 
			\left[ \begin{array}{c}
				\displaystyle \sum_{r=1}^R \langle \textbf{z}_1, \textbf{z}_r \rangle \cdot \textbf{x}_r \textbf{y}_1^T \cdot \textbf{v}_{Y_r}\\
				\vdots\\
				\displaystyle \sum_{r=1}^R \langle \textbf{z}_R, \textbf{z}_r \rangle \cdot \textbf{x}_r \textbf{y}_R^T \cdot \textbf{v}_{Y_r}
			\end{array} \right] = $$
			
			$$ = \left[ \begin{array}{c}
				\displaystyle \sum_{r=1}^R \langle \textbf{z}_1, \textbf{z}_r \rangle \textbf{x}_r \langle \textbf{y}_1, \textbf{v}_{Y_r} \rangle\\
				\vdots\\
				\displaystyle \sum_{r=1}^R \langle \textbf{z}_R, \textbf{z}_r \rangle \textbf{x}_r \langle \textbf{y}_R, \textbf{v}_{Y_r} \rangle
			\end{array} \right]$$
				
			$$ = 
			\left[ \begin{array}{c}
				\left[ \textbf{x}_1, \ldots, \textbf{x}_R \right]
				\left[ \begin{array}{c}
					\displaystyle \langle \textbf{z}_1, \textbf{z}_1 \rangle \langle \textbf{y}_1, \textbf{v}_{Y_1} \rangle\\
					\vdots\\
					\displaystyle \langle \textbf{z}_1, \textbf{z}_R \rangle \langle \textbf{y}_1, \textbf{v}_{Y_R} \rangle
				\end{array} 
				\right]\\
				\vdots\\
				\left[ \textbf{x}_1, \ldots, \textbf{x}_R \right]
				\left[ \begin{array}{c}
					\displaystyle \langle \textbf{z}_R, \textbf{z}_1 \rangle \langle \textbf{y}_R, \textbf{v}_{Y_1} \rangle\\
					\vdots\\
					\displaystyle \langle \textbf{z}_R, \textbf{z}_R \rangle \langle \textbf{y}_R, \textbf{v}_{Y_R} \rangle
				\end{array} 
				\right]
			\end{array} \right] = 
			\left[ \begin{array}{c}
				\textbf{X} \cdot \Big( (\Pi_{XY})_1 \ast \big( \textbf{V}_Y^T \cdot \textbf{y}_1 \big) \Big)\\
				\vdots\\
				\textbf{X} \cdot \Big( (\Pi_{XY})_R \ast \big( (\textbf{V}_Y^T \cdot \textbf{y}_R \big) \Big)
			\end{array} \right] = $$
			
			$$ = 
			vec\left( \left[ \textbf{X} \cdot \Big( (\Pi_{XY})_1 \ast \big( \textbf{V}_Y^T \cdot \textbf{y}_1 \big) \Big), \ldots, \textbf{X} \cdot \Big( (\Pi_{XY})_R \ast \big( (\textbf{V}_Y^T \cdot \textbf{y}_R \big) \Big) \right] \right) = $$
			
			$$ =
			vec\left( \textbf{X} \cdot \Big( \Pi_{XY} \ast \big( (\textbf{V}_Y^T \cdot \textbf{Y} \big) \Big) \right)$$
where each $(\Pi_{XY})_r$ is the $r$-th column of $\Pi_{XY}$. Despite the notation refers to the rows of $\Pi_{XY}$, this is not a problem since this matrix is symmetric.\\

			\textbf{Case 2 (equal modes):} In this case we have 
			
			$$\frac{\partial f}{\partial \textbf{X}}^T \frac{\partial f}{\partial \textbf{X}} \cdot 	vec(\textbf{V}_X) = $$	
			
			$$ = 
			\left[ \begin{array}{ccc}
				\displaystyle \langle \textbf{y}_1, \textbf{y}_1 \rangle \langle \textbf{z}_1, \textbf{z}_1 \rangle \cdot \textbf{I}_m & \ldots & \displaystyle \langle \textbf{y}_1, \textbf{y}_R \rangle \langle \textbf{z}_1, \textbf{z}_R \rangle \cdot \textbf{I}_m\\
				\vdots & & \vdots\\
				\displaystyle \langle \textbf{y}_R, \textbf{y}_1 \rangle \langle \textbf{z}_R, \textbf{z}_1 \rangle \cdot \textbf{I}_m & \ldots & \displaystyle \langle \textbf{y}_R, \textbf{y}_R \rangle \langle \textbf{z}_R, \textbf{z}_R \rangle \cdot \textbf{I}_m
			\end{array} \right]
			\left[ \begin{array}{c}
				\textbf{v}_{X_1}\\
				\vdots\\
				\textbf{v}_{X_R}
			\end{array}	\right] = $$

			$$ = \left[ \begin{array}{ccc}
				\displaystyle \sum_{r=1}^R \langle \textbf{y}_1, \textbf{y}_r \rangle \langle \textbf{z}_1, \textbf{z}_r \rangle \cdot \textbf{v}_{X_r}\\
				\vdots\\
				\displaystyle \sum_{r=1}^R \langle \textbf{y}_R, \textbf{y}_r \rangle \langle \textbf{z}_R, \textbf{z}_r \rangle \cdot \textbf{v}_{X_r}
			\end{array} \right] = $$
			
			$$\hspace{1cm} = 
			\left[ \begin{array}{c}
				\textbf{V}^{(\ell')} \cdot (\Pi_X)_1\\ 
				\vdots\\ 
				\textbf{V}^{(\ell')} \cdot (\Pi_X)_R
			\end{array} \right] = $$
			
			$$\hspace{1.9cm} = vec\left( \textbf{V}_X \cdot (\Pi_X)_1, \ldots, \textbf{V}^{(\ell')} \cdot (\Pi_X)_R \right) = vec\left( \textbf{V}_X \cdot \Pi_X \right). \hspace{1.9cm}\square$$

\section{Computational experiments}
\label{sec:comp}

\subsection{Procedure}
		We have selected a set of very distinct tensors to test the known tensor implementations. Given a tensor $\mathcal{T}$ and a rank $R$, we compute the CPD of TFX with the default maximum number of iterations\footnote{The default is $\texttt{maxiter} = 200$ iterations.} $100$ times and retain the best result, i.e., the CPD with the smallest relative error. Let $\varepsilon$ be this error. Now let ALG be any other algorithm implemented by some of the mentioned libraries. We set the maximum number of iterations to \verb|maxiter|, keep the other options with their defaults, and run ALG with these options $100$ times. The only accepted solutions are the ones with relative error smaller that $\varepsilon + \varepsilon/100$. Between all accepted solutions we return the one with the smallest running time. If none solution is accepted, we increase \verb|maxiter| by a certain amount and repeat. We try the values $\verb|maxiter|  = 5, 10, 50, 100, 150, \ldots, 900, 950, 1000$, until there is an accepted solution. The running time associated with the accepted solution is the accepted time. These procedures favour all algorithms against TFX since we are trying to find a solution close to the solution of TFX with the minimum number of iterations. We remark that the iteration process is always initiated with a random point. The option to generate a random initial point is offered by all libraries, and we use each one they offer (sometimes random initialization was already the default option). There is no much difference in their initializations, which basically amounts to draw random factor matrices from the standard Gaussian distribution. The time to perform the MLSVD or any kind of preprocessing is included in the time measurements. If one want to reproduce the tests presented here, they can be found at \url{https://github.com/felipebottega/Tensor-Fox/tree/master/tests}.
		
	\subsection{Algorithms}
		We used Linux Mint operational system in our tests. All tests were conducted using a processor Intel Core i7-4510U - 2.00GHz (2 physical cores and 4 threads) and 8GB of memory. The libraries mentioned run in Python or Matlab. We use Python - version 3.6.5 and Matlab - version 2017a. In both platforms we used BLAS MKL-11.3.1. Finally, we want to mention that TFX runs in Python using Numba to accelerate computations. We used the version 0.41 of Numba in these tests. The algorithm and implementations that are used in the experiments are the following.	
			
		\subsubsection{TFX}
			The algorithm used in TFX's implementation is the nonlinear squares scheme described in the previous section. There are more implementation details to be discussed, but the interested reader can check \cite{tensorfox} for more information. 
			
		\subsubsection{TALS}
			This is the Tensorlab's implementation of ALS algorithm. Although ALS is remarkably fast and easy to implement, it is not very accurate specially in the presence of bottlenecks or swamps. It seems (see \cite{tensorlab2}) that this implementation is very robust while still fast. 
				
		\subsubsection{TNLS} 
			This is the Tensorlab's implementation of NLS algorithm. This is the one we described in the previous section. We should remark that this implementation is similar to TFX's implementation at some points, but there are big differences when we look in more details. In particular the compression procedure, the preconditioner, the damping parameter update rule and the number of iterations of the conjugate gradient are very different.
				
		\subsubsection{TMINF}
			This is the Tensorlab's implementation of the problem as an optimization problem. They use a quasi-Newton method, the limited-memory BFGS, and consider equation ~\ref{eq:rank_approx_tucker} just as a minimization of a function. 
				
		\subsubsection{OPT}
			Just as the TMINF approach, the OPT algorithm is a implementation of Tensor Toolbox, which considers ~\ref{eq:rank_approx_tucker} as a minimization of a function. They claim that using the algorithm option 'lbfgs' is the preferred option\footnote{Check \url{https://www.tensortoolbox.org/cp_opt_doc.html} for more information.}, so we used this way.  
				
		\subsubsection{TLALS}
			TensorLy has only one way to compute the CPD, which is a implementation of the ALS algorithm. We denote it by TLALS, do not confuse with TALS, the latter is the Tensorlab's implementation.
			
		\subsubsection{fLMa}
			fLMA stands for \emph{fast Levenberg-Marquardt algorithm}, and it is a different version of the damped Gauss-Newton.\\
				
		For all Tensorlab algorithm implementation we recommend reading \cite{tensorlab2}, for the Tensor Toolbox we recommend \cite{cp_opt}, for TensorLy we recommend \cite{tensorly}, and for TensorBox we recommend \cite{tensorbox, cichoki2013}. In these benchmarks we used Tensorlab version 3.0 and Tensor Toolbox version 3.1.
				
		In Tensorlab it is possible to disable or not the option of refinement. The first action in Tensorlab is to compress the tensor and work with the compressed version. If we want to use refinement then the program uses the compressed solution to compute a refined solution in the original space. This can be more demanding but can improve the result considerably. In our experience working in the original space is not a good idea because the computational cost increases drastically and the gain in accuracy is generally very small. Still we tried all Tensorlab algorithms with and without refinement. We will write TALSR,  TNLSR and TMINFR for the algorithms TALS, TNLS and TMINF with refinement, respectively.  	
	
	\subsection{Tensors}		
		Now we describe the tensors used in our benchmarking. The idea was to have a set of very distinct tensors so we could test the implementations in very different situations. This can give a good idea of how they should perform in general. 
		
		\subsubsection{Swimmer} 
\begin{figure}[htbp] 
				\centering
				\includegraphics[scale=.4]{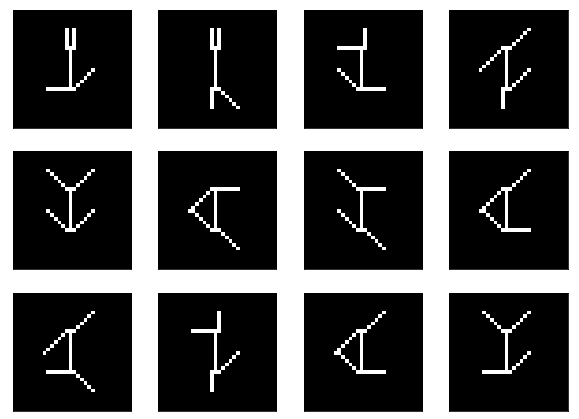}
				\caption{Swimmer tensor.}
				\label{fig:swimmer}
			\end{figure}

			This tensor was constructed based on the paper \cite{shashua} as an example of a non-negative tensor. It is a set of 256 images of dimensions $32 \times 32$ representing a swimmer. Each image contains a torso (the invariant part) of 12 pixels in the center and four limbs of 6 pixels that can be in one of 4 positions. In this work they proposed to use a rank $R = 50$ tensor to approximate, and we do the same for our test. In figure~\ref{fig:swimmer} we can see some frontal slices of this tensor.
			
		\subsubsection{Handwritten digit} 
\begin{figure}[htbp] 
				\centering
				\includegraphics[scale=.4]{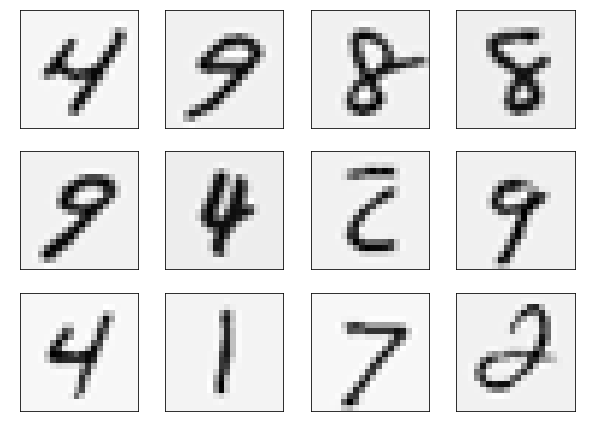}
				\caption{Handwritten digits tensor.}
				\label{fig:handwritten}
			\end{figure}

			This is a classic tensor in machine learning, it is the MNIST\footnote{\url{http://yann.lecun.com/exdb/mnist/}} database of handwritten digits. Each slice is a image of dimensions $20 \times 20$ of a handwritten digit. Also, each 500 consecutive slices correspond to the same digit, so the first 500 slices correspond to the digit 0, the slices 501 to 1000 correspond to the digit 1, and so on. We choose $R = 150$ as a good rank to construct the approximating CPD to this tensor. In figure~\ref{fig:handwritten} we can see some frontal slices of this tensor.
			
		\subsubsection{Border rank} 
			We say a tensor $\mathcal{T}$ has \emph{border} rank $\underline{R}$ if there is a sequence of tensors of rank $\underline{R}$ converging to $\mathcal{T}$ and also there is not a sequence of tensors with lower rank satisfying the same property. If $\mathcal{T}$ has rank equal to $R$, then it is easy to see that $\underline{R} \leq R$. The strict inequality can happen too, and this means that the set of rank $R$ tensors is not closed. This phenomenon makes the CPD computation a challenging problem. This subject was discussed by Silva and Lim in \cite{lim}. In the same paper they showed that 
			$$\mathcal{T}_k = k \left( \textbf{x}_1 + \frac{1}{k} \textbf{y}_1 \right) \otimes \left( \textbf{x}_2 + \frac{1}{k} \textbf{y}_2 \right) \otimes \left( \textbf{x}_3 + \frac{1}{k} \textbf{y}_3 \right) - k \textbf{x}_1 \otimes \textbf{x}_2 \otimes \textbf{x}_3$$
is a sequence of rank sequence of rank 2 tensors converging to a tensor of rank 3, where each pair $\textbf{x}_i, \textbf{y}_i \in \mathbb{R}^m$ is linearly independent. The limit tensor is $\mathcal{T} = \textbf{x}_1 \otimes \textbf{x}_2 \otimes \textbf{y}_3 + \textbf{x}_1 \otimes \textbf{y}_2 \otimes \textbf{x}_3 + \textbf{y}_1 \otimes \textbf{x}_2 \otimes \textbf{x}_3$. We choose to compute a CPD of rank $R = 2$ to see how the algorithms behaves when we try to approximate a problematic tensor by tensor with low rank. In theory it is possible to have arbitrarily good approximations.

		\subsubsection{Matrix multiplication}
			Let $\mathcal{M}_N \in \mathbb{R}^{N^2 \times N^2 \times N^2}$ be the tensor associated with the multiplication between two matrices in $\mathbb{R}^{N \times N}$. The classic form of $\mathcal{M}_N$ is given by
			$$\mathcal{M}_N = \sum_{i=1}^m \sum_{j=1}^n \sum_{k=1}^l vec(\textbf{e}_j^i) \otimes vec(\textbf{e}_k^j) \otimes vec(\textbf{e}_k^i),$$
where $\textbf{e}_i^j$ is the matrix $N \times N$ with entry $(i,j)$ equal to 1 and the remaining entries equal to zero. Since Strassen \cite{strassen} it is known that matrix multiplication between matrices of dimensions $N \times N$ can be made with $\mathcal{O}(N^{\log_2 7})$ operations. Many improvements were made after Strassen but we won't enter in the details here. For the purpose of testing we choose the small value $N = 5$ and the rank $R = \lceil 5^{\log_2 7} \rceil = 92$. This value is the bound obtained by Strassen in \cite{strassen}. It is not necessarily the true rank of the tensor but it is close enough to make an interesting test.   

		\subsection{Collinear factors}
			The phenomenon of swamps occurs when all factors in each mode are almost collinear. Their presence is a challenge for many algorithms because they can slow down convergence. Now we will create synthetic data to simulate various degrees of collinearity between the factors. We begin generating three random matrices $\textbf{M}_X \in \mathbb{R}^{m \times R}, \textbf{M}_Y \in \mathbb{R}^{n \times R}, \textbf{M}_Z \in \mathbb{R}^{p \times r}$, where each entry is drawn from the normal distribution with mean 0 and variance 1. After that we perform QR decomposition of each matrix, obtaining the decompositions $\textbf{M}_X = \textbf{Q}_X \textbf{R}_X, \textbf{M}_Y = \textbf{Q}_Y \textbf{R}_Y, \textbf{M}_Z = \textbf{Q}_Z \textbf{R}_Z$. The matrices $\textbf{Q}_X, \textbf{Q}_Y, \textbf{Q}_Z$ are orthogonal. Now fix three columns $\textbf{q}_X^{(i')}, \textbf{q}_Y^{(j')}, \textbf{q}_Z^{(k')}$ of each one of these matrices. The factors $\textbf{X} = [ \textbf{x}^{(1)}, \ldots, \textbf{x}^{(R)} ] \in \mathbb{R}^{m \times R}$, $\textbf{Y} = [ \textbf{y}^{(1)}, \ldots, \textbf{y}^{(R)} ] \in \mathbb{R}^{n \times R}$, $\textbf{Z} = [ \textbf{z}^{(1)}, \ldots, \textbf{z}^{(R)} ] \in \mathbb{R}^{p \times R}$ are generated by the equations below.
			
			$$\textbf{x}^{(i)} = \textbf{q}_X^{(i')} + c \cdot \textbf{q}_X^{(i)}, \quad i=1 \ldots R$$
			$$\textbf{y}^{(j)} = \textbf{q}_Y^{(j')} + c \cdot \textbf{q}_Y^{(j)}, \quad j=1 \ldots R$$
			$$\textbf{z}^{(k)} = \textbf{q}_Z^{(k')} + c \cdot \textbf{q}_Z^{(k)}, \quad k=1 \ldots R$$  
		
			The parameter $c \geq 0$ defines the degree of collinearity between the vectors of each factor. A value of $c$ close to 0 indicates high degree of collinearity, while a high value of $c$ indicates low degree of collinearity.
			
			Another phenomenon that occurs in practice is the presence of noise in the data. So we will treat these two phenomena at once in this benchmark. After generating the factors $\textbf{X}, \textbf{Y}, \textbf{Z}$ we have a tensor $\mathcal{T} = (\textbf{X}, \textbf{Y}, \textbf{Z}) \cdot \mathcal{I}_{R \times R \times R}$. That is, $\textbf{X}, \textbf{Y}, \textbf{Z} $ are the exact CPD of $\mathcal{T}$. Now consider a noise $\mathcal{N} \in \mathbb{R}^{m \times n \times p}$ such that each entry of $\mathcal{N} $ is obtained by the normal distribution with mean 0 and variance 1. Thus we form the tensor $\hat{\mathcal{T}} = \mathcal{T} + \nu \cdot \mathcal{N}$, where $ \nu > 0 $ defines the magnitude of the noises. The idea is to compute a CPD of $ \hat{\mathcal{T}}$ of rank $R$ and then evaluate the relative error between this tensor and $\mathcal{T}$. We expect the computed CPD to clear the noises and to be close to $\mathcal{T}$ (even if it is not close to $ \hat{\mathcal{T}}$). We will fix $\nu = 0.01$ and generate tensors for $c = 0.1, 0.5, 0.9$. In all cases we will be using $ m = n = p = 300 $ and $R = 15 $. This is a particularly difficult problem since we are considering swamps and noises at once. The same procedure to generate tensors were used for benchmarking in \cite{cichoki2013}.
			
		\subsection{Double bottlenecks}
			We proceed almost in the same as before for swamps, we used the same procedure to generate the first two columns of each factor matrix, then the remaining columns are equal to the columns of the QR decomposition. After generating the factors $\textbf{X}, \textbf{Y}, \textbf{Z}$ we consider a noise $\mathcal{N} \in \mathbb{R}^{m \times n \times p}$ such that each entry of $\mathcal{N} $ is obtained by the normal distribution with mean 0 and variance 1. Thus we form the tensor $\hat{\mathcal{T}} = \mathcal{T} + \nu \cdot \mathcal{N}$, where $ \nu > 0 $ defines the magnitude of the noises. The collinear parameter used for the tests are $c = 0.1, 0.5$, and we fix $\nu = 0.01$ as before. The procedure before the noise is presented below.
			
			$$\textbf{x}^{(i)} = \textbf{q}_X^{(i')} + c \cdot \textbf{q}_X^{(i)}, \quad i=1, 2$$
			$$\textbf{y}^{(j)} = \textbf{q}_Y^{(j')} + c \cdot \textbf{q}_Y^{(j)}, \quad j=1, 2$$
			$$\textbf{z}^{(k)} = \textbf{q}_Z^{(k')} + c \cdot \textbf{q}_Z^{(k)}, \quad k=1, 2$$ 
			
			$$\textbf{x}^{(i)} = \textbf{q}_X^{(i)}, \quad i=3 \ldots R$$
			$$\textbf{y}^{(j)} = \textbf{q}_Y^{(j)}, \quad j=3 \ldots R$$
			$$\textbf{z}^{(k)} = \textbf{q}_Z^{(k)}, \quad k=3 \ldots R$$ 
			
	\subsection{Results}		
\begin{figure}[htbp] 
			\includegraphics[scale=0.3]{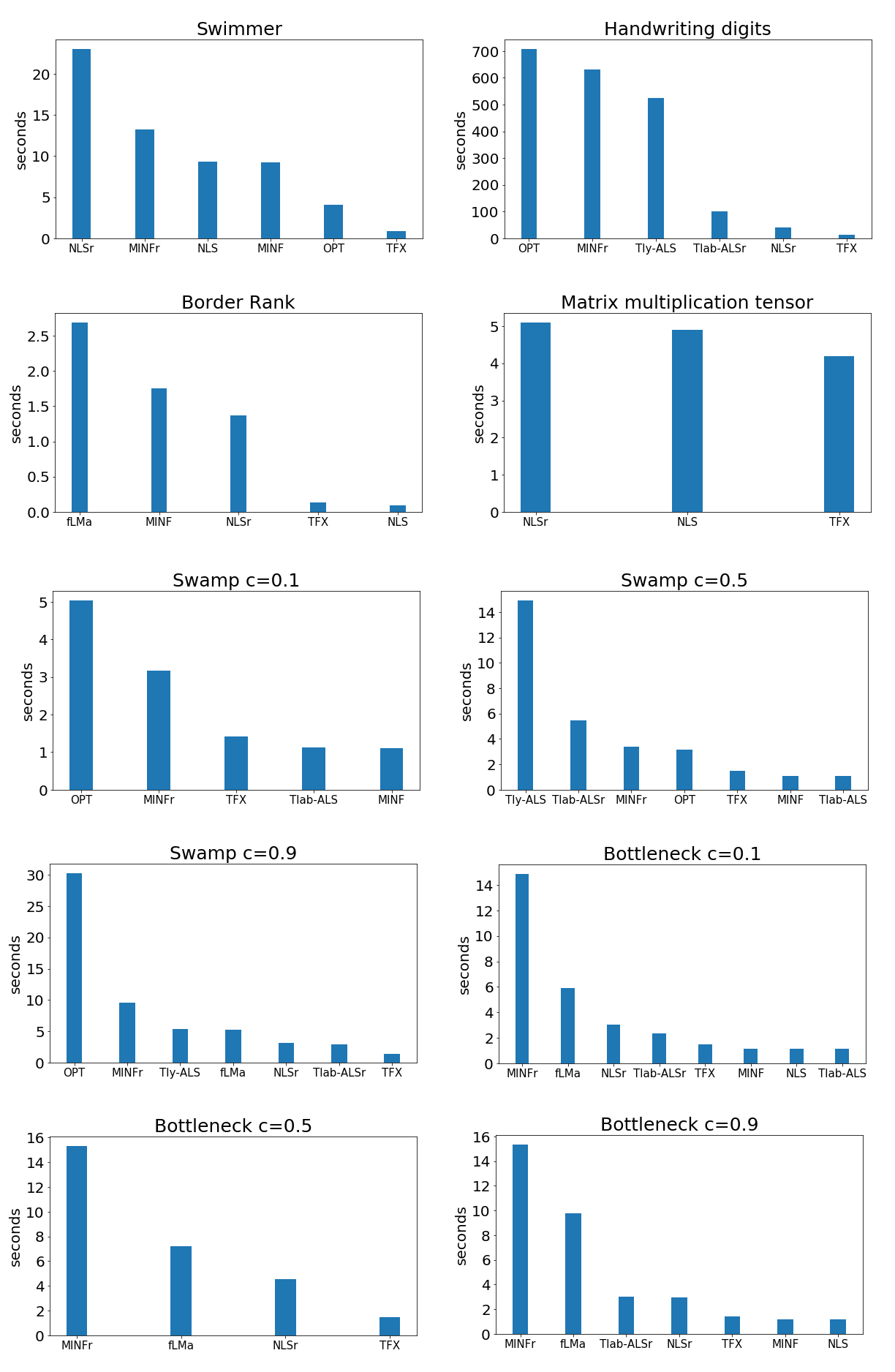}
			\caption{Benchmarks of all tensors and all implementations.}
			\label{fig:benchs}
		\end{figure}

		In figure~\ref{fig:benchs} there are some charts, each one showing the best running time of the algorithms with respect to each one of the tensors describe previously. If some algorithm is not included in a chart, it means that the algorithm was unable to achieve an acceptable error within the conditions described at the beginning if this section. 
		
		The first thing we should note is that TNLS is robust enough so it can deliver an acceptable solution in almost every test, with the exception of the border rank tensor and the collinear tensor with $c = 0.9$. Not only that but it is also very fast, beating TFX in two tests. It should be noted too how TALS performs well at the collinear tests. Normally we could expect the opposite since there are many swamps at these tests. This shows how Tensorlab made a good work with the ALS algorithm. Apart from that, none of the other algorithms seemed to stand out for some test. Finally, we want to add that, although none algorithm was able to deliver an acceptable solution for the border rank tensor test, the algorithm TNLS could compute solutions very close to the desired interval within a reasonable time (but still more time than TFX).

\section{Gaussian Mixture}
\label{sec:gmix}

Gaussian mixture models are more connected to real applications. We start briefly describing the model, them we apply some algorithms in synthetic data and compare the results. The theory discussed here is based on \cite{anandkumar}.
	
	\subsection{Gaussian mixture model} 
		Consider a mixture of $K$ Gaussian distributions with identical covariance matrices. We have lots of data with unknown averages and unknown covariance matrices to infer the parameters of distributions, which are are unknown. The problem at hand is to design an algorithm to \emph{learn} these parameters from the data given. We use $\mathbb{P}$ to denote probability and $\mathbb{E}$ to denote expectation (which we may also call the \emph{mean} or \emph{average}).
		
\begin{figure}[htbp] 
			\centering
			\includegraphics[scale=.42]{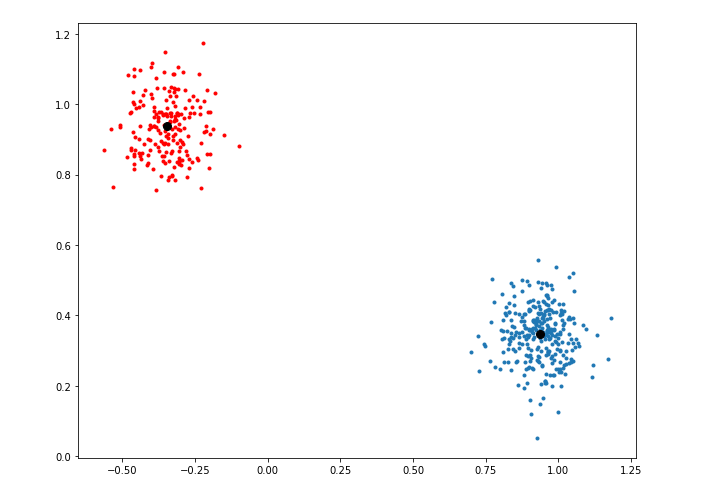}
			\caption{Gaussian mixture in the plane with 2 distributions. The first distribution has mean $\textbf{u}^{(1)} = [-0.34, \  0.93]^T$ and the second has mean $\textbf{u}^{(2)} = [0.93, \ 0.34]^T$. The variance is $\sigma^2 = 0.0059$.}
			\label{fig:gaussian_mix}
		\end{figure}

		Let $\textbf{x}^{(1)}, \ldots, \textbf{x}^{(N)} \in \mathbb{R}^d$ be a set of collected data sample. Let $h$ be a discrete random variable with values in $\{1, 2, \ldots, K\}$ such that $\mathbb{P}[h = i]$ is the probability that a sample $\textbf{x}$ is a member of the $i$-th distribution. We denote $w^{(i)} = \mathbb{P}[h = i]$ and $\textbf{w} = [ w^{(1)}, \ldots, w^{(K)}]^T$, the vector of probabilities. Let $\textbf{u}^{(i)} \in \mathbb{R}^d$ be the mean of the $i$-th distribution and assume that all distributions have the same covariance matrix $\sigma^2 \textbf{I}_d$ for $\sigma > 0$. See figure~\ref{fig:gaussian_mix} for an illustration of a Gaussian mixture in the case where $d = 2$ and $K = 2$.
		
		To keep everything simple we also assume that the means $\textbf{u}^{(i)}$ form an orthonormal set, i.e., every $\textbf{u}^{(i)}$ is a unit vector and these vectors are orthogonal with respect to each other. To see how to proceed in the general case we refer to the previously mentioned paper. 
		
		Given a sample point $\textbf{x}$, note that we can write
		$$\textbf{x} = \textbf{u}_h + \textbf{z},$$
where $\textbf{z}$ is a random vector with mean 0 and covariance $\sigma^2 \textbf{I}_d$. We summarise the main results in the next theorem whose proof can be found in \cite{anandkumar}.

		\begin{theorem}[Hsu and Kakade, 2013] \label{moments} 
			Assume $d \geq K$. The variance $\sigma^2$ is the smallest eigenvalue of the covariance matrix $\mathbb{E}[ \textbf{x} \otimes \textbf{x} ] - \mathbb{E}[ \textbf{x} ] \otimes \mathbb{E}[ \textbf{x} ]$. Furthermore, if		
			\begin{flalign*}
				& M_1 = \mathbb{E}[ \textbf{x} ],\\
				& M_2 = \mathbb{E}[ \textbf{x} \otimes \textbf{x} ] - \sigma^2 \textbf{I}_d,\\
				& M_3 = \mathbb{E}[ \textbf{x} \otimes \textbf{x} \otimes \textbf{x} ] - \sigma^2 \sum_{i=1}^d \left( \mathbb{E}[ \textbf{x} ] \otimes \textbf{e}_i \otimes \textbf{e}_i + \textbf{e}_i \otimes \mathbb{E}[ \textbf{x} ] \otimes \textbf{e}_i + \textbf{e}_i \otimes \textbf{e}_i \otimes \mathbb{E}[ \textbf{x} ] \right),
			\end{flalign*}
then
			\begin{flalign*}
				& M_1 = \sum_{i=1}^K w^{(i)}\ \textbf{u}^{(i)},\\
				& M_2 = \sum_{i=1}^K w^{(i)}\ \textbf{u}^{(i)} \otimes \textbf{u}^{(i)},\\
				& M_3 = \sum_{i=1}^K w^{(i)}\ \textbf{u}^{(i)} \otimes \textbf{u}^{(i)} \otimes \textbf{u}^{(i)}.
			\end{flalign*}
		\end{theorem} 
		
		Theorem~\ref{moments} allows us to use the method of moments, which is a classical parameter estimation technique from statistics. This method consists in computing certain statistics of the data (often empirical moments) and use it to find model parameters that give rise to (nearly) the same corresponding population quantities. Now suppose that $N$ is large enough so we have a reasonable number of sample points to make useful statistics. First we compute the empirical mean 
		
		\begin{equation}		
			\hat{\mu} := \frac{1}{N} \sum_{j=1}^N \textbf{x}^{(j)} \approx \mathbb{E}[ \textbf{x} ]. \label{eq:empirical_mean}
		\end{equation}
		
		Now use this result to compute the empirical covariance matrix
		
		\begin{equation}
			\hat{\textbf{S}} := \frac{1}{N} \sum_{j=1}^N ( \textbf{x}^{(j)} \otimes \textbf{x}^{(j)} - \hat{\mu} \otimes \hat{\mu} ) \approx \mathbb{E}[ \textbf{x} \otimes \textbf{x} ] - \mathbb{E}[ \textbf{x} ] \otimes \mathbb{E}[ \textbf{x} ]. \label{eq:empirical_cov}
		\end{equation}
		
		The smallest eigenvalue of $\hat{\textbf{S}}$ is the empirical variance $\hat{\sigma}^2 \approx \sigma^2$. Now we compute the empirical third moment (empirical skewness)
		
		 \begin{equation}
			\hat{\mathcal{S}} := \frac{1}{N} \sum_{j=1}^N \textbf{x}^{(j)} \otimes \textbf{x}^{(j)} \otimes \textbf{x}^{(j)} \approx \mathbb{E}[ \textbf{x} \otimes \textbf{x} \otimes \textbf{x} ] \label{eq:empirical_skew}
		\end{equation}
and use it to get the empirical value of $M_3$,
		\begin{equation}
			\hat{\mathcal{M}}_3 := \hat{\mathcal{S}} - \hat{\sigma}^2 \sum_{i=1}^d \left( \hat{\mu} \otimes \textbf{e}_i \otimes \textbf{e}_i + \textbf{e}_i \otimes \hat{\mu} \otimes \textbf{e}_i + \textbf{e}_i \otimes \textbf{e}_i \otimes \hat{\mu} \right) \approx M_3.\label{eq:empirical_M3}
		\end{equation}
		
		By theorem~\ref{moments}, $M_3 = \displaystyle \sum_{i=1}^K w^{(i)}\ \textbf{u}^{(i)} \otimes \textbf{u}^{(i)} \otimes \textbf{u}^{(i)}$, which is a symmetric tensor containing all parameter information we want to find. The idea is, after computing a symmetric CPD for $\hat{\mathcal{M}}_3$, normalize the factors so each vector has unit norm. By doing this we have a tensor of the form
		$$\sum_{i=1}^K \hat{w}^{(i)}\ \hat{\textbf{u}}^{(i)} \otimes \hat{\textbf{u}}^{(i)} \otimes \hat{\textbf{u}}^{(i)}$$
as a candidate to solution. Note that it is easy to make all $\hat{w}^{(i)}$ positive. If some of them is negative, just multiply it by $-1$ and multiply one of the associated vectors also by $-1$. The final tensor is unchanged but all $\hat{w}^{(i)}$ now are positive.   

	\subsection{Computational experiments}
		Now we describe how to generated data, what algorithms are used to compute the CPDs and how the comparisons are made. Here we restrict our attention only to the implementations of TFX and Tensorlab. 
		
		To compute a symmetric CPD with TFX we just have to set the option \verb|symm| to\footnote{\url{https://github.com/felipebottega/Tensor-Fox/blob/master/tutorial/3-intermediate_options.ipynb}} \verb|True|. We also observed that it was necessary to let TFX to perform more conjugate gradient iterations in order to obtain meaningful results, although this may increase the computational time. To compute a symmetric CPD with Tensorlab we need to create a model specifying that the result should be a symmetric tensor. Also, in this context it is only possible to use NLS and MINF algorithms, with or without refinement. For more information we recommend reading sections 8.2 and 8.3 from the Tensorlab guide.\footnote{\url{https://www.tensorlab.net/userguide3.pdf}}	In order to obtain meaningful results	we set the parameters \verb|TolFun| and \verb|TolX| to $10^{-12}$. 
		
		\subsubsection{Procedure}
			To work with Gaussian mixtures we generate datasets for $d = 20, 100$ and $K = 5, 15$. In particular, we have that $\mathcal{M}_3 \in \mathbb{R}^{d \times d \times d}$. For each example we generate the probabilities $w^{(i)}$ by taking random values in the interval $(0,1)$ such that all $w^{(i)} > 0$ and $w^{(1)} + \ldots + w^{(K)} = 1$. To generate the means we first generated a random Gaussian matrix $\textbf{M} \in \mathbb{R}^{d \times K}$ with full rank, computed its SVD, $\textbf{M} = \textbf{U} \Sigma \textbf{V}^T$, and used each column of $\textbf{U}$ to be a mean of the distribution (we only need the first $K$ columns of $\textbf{U}$).  
			
			For each example we generated a population with size $N = 10000$ by drawing samples $\textbf{x} \in \mathbb{R}^d$ using the following procedure:
			\begin{enumerate}
				\item Generate a random number $h \in \{1, 2, \ldots, K\}$ from the distribution given by the $w^{(i)}$. More precisely, $\mathbb{P}[h = i] = w^{(i)}$. The number obtained, $i$, is the distribution of $\textbf{x}$.
				\item Generate a random vector $\textbf{z} \in \mathbb{R}^d$ with mean 0 and covariance $\sigma^2 \textbf{I}_d$.
				\item Set $\textbf{x} = \textbf{u}^{(i)} + \textbf{z}$.
				\item Repeat the previous steps $N$ times.
			\end{enumerate}	
			 
			With these samples we are able to estimate $\mathcal{M}_3$ using the empirical tensor $\hat{\mathcal{M}}_3$ defined in~\ref{eq:empirical_M3}.
					 
		\subsubsection{Computational results}
			For each example we computed 100 CPDs and retain the one with best fit. In this case the best fit is defined as the smallest CPD error but the error corresponding to the parameters. If $\hat{w}^{(i)}$ and $\hat{\textbf{u}}^{(i)}, i=1 \ldots K$, are the approximated parameters, then the corresponding fit is the value 
			$$\frac{ \| \hat{\textbf{w}} -  \textbf{w} \| }{ \| \textbf{w} \| } + \frac{ \| \hat{\textbf{U}} - \textbf{U} \| }{ \| \textbf{U} \| },$$
where $\hat{\textbf{w}} = [ \hat{w}^{(1)}, \ldots, \hat{w}^{(K)}]^T$ and $\hat{\textbf{U}} = [ \hat{\textbf{u}}^{(1)}, \ldots, \hat{\textbf{u}}^{(K)} ]$.

			For this problem we took a different approach when comparing the algorithms. We let all programs run with the parameters mentioned before and compared the results in a \verb|accuracy| $\times$ \verb|time| plot. For each example we make two plots, the first with the errors of the probabilities and the second with the errors of the means.		
			
\begin{figure}[htbp]
				\includegraphics[scale=0.4]{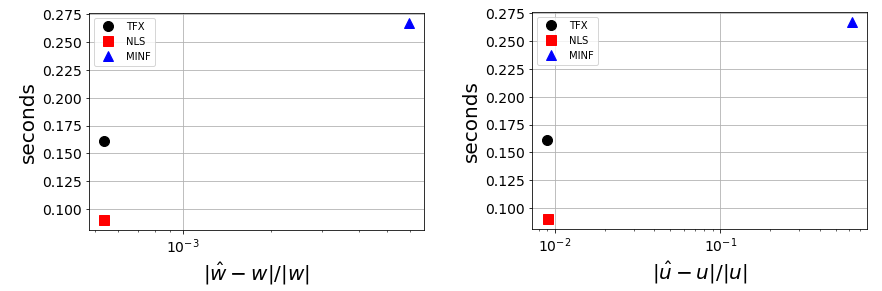}
				\caption{Computational results for $d = 20$ and $K = 5$. The horizontal axis represents the average time (in seconds) to compute the solutions. At the left, the vertical axis represents the best relative error of the probabilities, and at the right, the vertical axis represents the best relative error of the means.}
				\label{fig:d20_k5}
			\end{figure}

\begin{figure}[htbp]
				\includegraphics[scale=0.415]{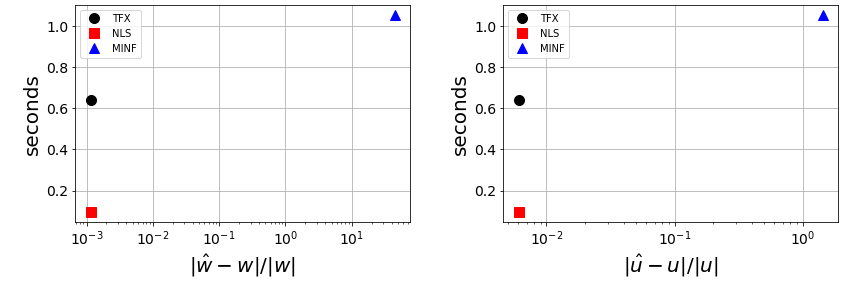}
				\caption{Computational results for $d = 20$ and $K = 15$.} 
				\label{fig:d20_k15}
			\end{figure}

\begin{figure}[htbp]
				\includegraphics[scale=0.408]{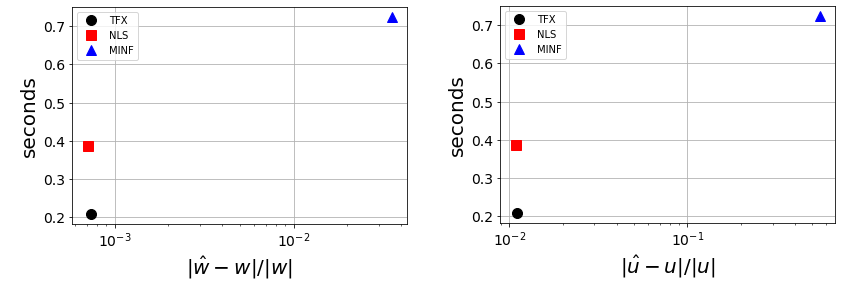}
				\caption{Computational results for $d = 100$ and $K = 5$.}
				\label{fig:d100_k5}
			\end{figure}

\begin{figure}[htbp]
				\includegraphics[scale=0.4]{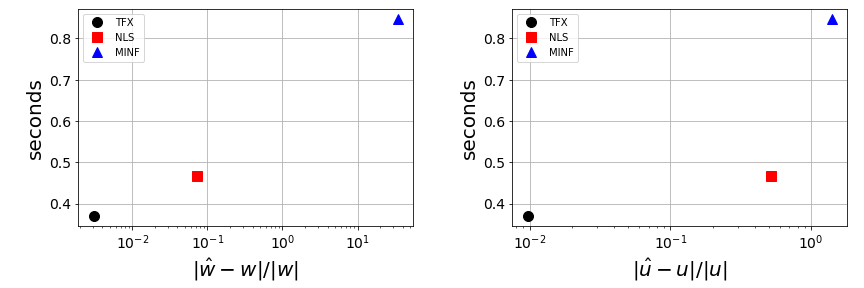}
				\caption{Computational results for $d = 100$ and $K = 15$.}
				\label{fig:d100_k15}
			\end{figure}

			For $d = 20$ we can see that both NLS and TFX achieves the same accuracy, with NLS being a little faster. When the dimension is increased to $d = 100$, TFX starts to be faster than NLS. Not only that, but for $K = 15$ we can see that NLS is less accurate than TFX. The difference in speed and accuracy becomes more evident as $d$ and $K$ increases. In all tests MINF performed poorly, being the slowest and less accurate. See figures~\ref{fig:d20_k5},~\ref{fig:d20_k15},~\ref{fig:d100_k5},~\ref{fig:d100_k15}.
			
\section{Conclusions}
\label{sec:conclusions}

		In the past years several implementations of algorithms to compute the CPD were proposed. Usually they are based on alternating least squares, nonlinear squares or optimization methods. Our work differs from others due to several other implementation features not deeply investigated so far: the damping parameter rule update, preconditioners, regularization, compression, strategies for the number of conjugate gradient iterations, and others. In order to construct TensorFox we conducted a deep investigation through several possibilities at these features. After several attempts and tests we converged to this competitive algorithm for computing the CPD. 
		
		We introduced some of the main features of our implementation and performed a series of tests of it against other implementations. We then introduced the concept of Gaussian mixtures and performed more tests. This particular problem is harder than the others since we need to compute several CPDs in order have a good fit. In these tests our algorithm also showed to be competitive and, in particular, it seemed to perform better as the problem gets harder. 

\appendix
\section{Proofs and some generalization} 
Although theorems 4.3 and 4.4 were stated for third order tensors, their generalizations are interesting and can be easily addressed. We remark that theorem 4.6 also can be generalized without much effort but its third order version is the one of interest in this paper.
	
	Denote $\textbf{w} = [vec(\textbf{W}^{(1)})^T, \ldots, vec(\textbf{W}^{(L)})^T]^T$, where each $\textbf{W}^{(\ell)} \in \mathbb{R}^{I_\ell \times R}$ is the factor matrix of a $L$-th order CPD given by $(\textbf{W}^{(1)}, \ldots, \textbf{W}^{(L)}) \cdot \mathcal{I}_{R \times \ldots \times R}$.\bigskip

	\textbf{Proof of theorem 4.3:} To prove item 1, just take any $\textbf{w} \in \mathbb{R}^{R \sum_{\ell=1}^L I_\ell}$ and note that
		
		$$\left\langle \left( \textbf{J}_f^T \textbf{J}_f + \mu \textbf{D} \right)\textbf{w}, \textbf{w} \right\rangle = $$
		$$ = \left\langle \textbf{J}_f^T \textbf{J}_f \textbf{w}, \textbf{w} \right\rangle + \left\langle \mu \textbf{D} \textbf{w}, \textbf{w} \right\rangle = $$
		$$ = \left\langle \textbf{J}_f \textbf{w}, \textbf{J}_f \textbf{w} \right\rangle + \left\langle \sqrt{\mu} \sqrt{\textbf{D}} \textbf{w}, \sqrt{\mu} \sqrt{\textbf{D}} \textbf{w} \right\rangle = $$
		$$ = \| \textbf{J}_f \textbf{w} \|^2 + \| \sqrt{\mu} \sqrt{\textbf{D}} \textbf{w} \|^2 > 0.$$
	
		The proof of item 2 is very similar to the previous proof in the classic Gauss-Newton. From the iteration formula
		
		$$\textbf{w}^{(k+1)} = \textbf{w}^{(k)} - \left( \textbf{J}_f(\textbf{w}^{(k)})^T \textbf{J}_f(\textbf{w}^{(k)}) + \mu \textbf{D} \right)^{-1} \textbf{J}_f(\textbf{w}^{(k)})^T \cdot f(\textbf{w}^{(k)})$$
we can conclude that

		$$- \left( \textbf{J}_f(\textbf{w}^{(k)})^T \textbf{J}_f(\textbf{w}^{(k)}) + \mu \textbf{D} \right) \cdot \left( \textbf{w}^{(k+1)} - \textbf{w}^{(k)} \right) = \textbf{J}_f(\textbf{w}^{(k)})^T \cdot f(\textbf{w}^{(k)}).$$
		
		Now, with this identity, note that
		
		$$\left\langle \nabla F(\textbf{w}^{(k)}),  \textbf{w}^{(k+1)} - \textbf{w}^{(k)} \right\rangle = $$
		$$ = \left\langle \textbf{J}_f(\textbf{w}^{(k)})^T f(\textbf{w}^{(k)}), \textbf{w}^{(k+1)} - \textbf{w}^{(k)} \right\rangle = $$
		$$ = - \left\langle \left( \textbf{J}_f(\textbf{w}^{(k)})^T \textbf{J}_f(\textbf{w}^{(k)}) + \mu \textbf{D} \right) \cdot \left( \textbf{w}^{(k+1)} - \textbf{w}^{(k)} \right), \textbf{w}^{(k+1)} - \textbf{w}^{(k)} \right\rangle < 0. $$
		
		The inequality above follows from the fact that $\textbf{J}_f(\textbf{w}^{(k)})^T \textbf{J}_f(\textbf{w}^{(k)}) + \mu \textbf{D}$ is positive definite.
	
		To prove item 3, take $\mu$ such that $\| \textbf{D}^{-1} \textbf{J}_f(\textbf{w}^{(k)})^T \textbf{J}_f(\textbf{w}^{(k)}) \| \ll \mu$ (this is ``large enough'' in this context). We know $\textbf{J}_f(\textbf{w}^{(k)})^T \textbf{J}_f(\textbf{w}^{(k)}) + \mu \textbf{D}$ since it is positive definite. Also, by the definition of $\mu$ we have that
		
		$$\left( \textbf{J}_f(\textbf{w}^{(k)})^T \textbf{J}_f(\textbf{w}^{(k)}) + \mu \textbf{D} \right)^{-1} = \left( \mu \textbf{D} \left( \displaystyle \frac{1}{\mu} \textbf{D}^{-1} \textbf{J}_f(\textbf{w}^{(k)})^T \textbf{J}_f(\textbf{w}^{(k)}) + \textbf{I} \right) \right)^{-1} \approx $$
		$$ \approx \left( \mu \textbf{D} \left( \textbf{0} + \textbf{I} \right) \right)^{-1} = \frac{1}{\mu} \textbf{D}^{-1}.$$
		
		Using the iteration formula with this approximation gives 
		$$\textbf{w}^{(k+1)} \approx \textbf{w}^{(k)} - \frac{1}{\mu} \textbf{D}^{-1} \textbf{J}_f(\textbf{w}^{(k)})^T f(\textbf{w}^{(k)}) = \textbf{w}^{(k)} - \frac{1}{\mu} \textbf{D}^{-1} \nabla F(\textbf{w}^{-1}).$$	

		Finally, to prove item 4 just consider $\mu \approx 0$ and substitute in the iteration formula. Then we get the classical formula trivially.$\hspace{6.2cm}\square$\bigskip
		
	Let two matrices $\textbf{A} \in \mathbb{R}^{k \times \ell}, \textbf{B} \in \mathbb{K}^{m \times n}$. The \emph{Kronecker product} between $\textbf{A}$ and $\textbf{B}$ is defined by
				$$\textbf{A} \tilde{\otimes} \textbf{B} = \left[
				\begin{array}{cccc}
					a_{11} \textbf{B} & a_{12} \textbf{B} & \ldots & a_{1\ell} \textbf{B}\\
					a_{21} \textbf{B} & a_{22} \textbf{B} & \ldots & a_{2\ell} \textbf{B}\\
					\vdots & \vdots & \ddots & \vdots\\
					a_{k1} \textbf{B} & a_{k2} \textbf{B} & \ldots & a_{k\ell} \textbf{B}
				\end{array}
				\right].$$
	
			The matrix given in the definition is a block matrix such that each block is a $m \times n$ matrix, so $\textbf{A} \tilde{\otimes} \textbf{B}$ is a $km \times \ell n$ matrix. We would like to point out that some texts uses $\otimes$ for the Kronecker product and $\circ$ for the tensor product.\bigskip 
		
	\textbf{Theorem 4.4 generalized:}
				Denote $\omega_{r' r''}^{(\ell)} = \langle \textbf{w}_{r'}^{(\ell)}, \textbf{w}_{r''}^{(\ell)} \rangle$. Then we have that			
			
				$$\textbf{J}_f^T \textbf{J}_f = 
					\left[ \begin{array}{ccc}
						\textbf{H}_{11} & \ldots & \textbf{H}_{1L}\\
						\vdots & & \vdots\\
						\textbf{H}_{L1} & \ldots & \textbf{H}_{LL}
					\end{array} \right],$$
where

				$$\textbf{H}_{\ell' \ell''} = 
					\left[ \begin{array}{ccc}
					\displaystyle \prod_{\ell \neq \ell', \ell''} \omega_{11}^{(\ell)} \cdot \textbf{w}_1^{(\ell')} \big( \textbf{w}_1^{(\ell'')} \big)^T & \ldots & \displaystyle \prod_{\ell \neq \ell', \ell''} \omega_{1R}^{(\ell)} \cdot \textbf{w}_R^{(\ell')} \big( \textbf{w}_1^{(\ell'')} \big)^T\\
					\vdots & & \vdots\\
					\displaystyle \prod_{\ell \neq \ell', \ell''} \omega_{R1}^{(\ell)} \cdot \textbf{w}_1^{(\ell')} \big( \textbf{w}_R^{(\ell'')} \big)^T & \ldots & \displaystyle \prod_{\ell \neq \ell', \ell''} \omega_{RR}^{(\ell)} \cdot \textbf{w}_R^{(\ell')} \big( \textbf{w}_R^{(\ell'')} \big)^T
				\end{array} \right]$$
for $\ell' \neq \ell''$, and

				$$\textbf{H}_{\ell' \ell'} = 
					\left[ \begin{array}{ccc}
					\displaystyle \prod_{\ell \neq \ell'} \omega_{11}^{(\ell)} \cdot \textbf{I}_{I_{\ell'}} & \ldots & \displaystyle \prod_{\ell \neq \ell'} \omega_{1R}^{(\ell)} \cdot \textbf{I}_{I_{\ell'}}\\
					\vdots & & \vdots\\
					\displaystyle \prod_{\ell \neq \ell'} \omega_{R1}^{(\ell)} \cdot \textbf{I}_{I_{\ell'}} & \ldots & \displaystyle \prod_{\ell \neq \ell'} \omega_{RR}^{(\ell)} \cdot \textbf{I}_{I_{\ell'}}
				\end{array} \right].$$\bigskip
				
			\textbf{Proof:} First note that
			
				$$\textbf{J}_f^T \textbf{J}_f = 
			 	\left[ \begin{array}{c}
			 		\displaystyle\frac{\partial f}{\partial \textbf{W}^{(1)}}^T\\ 
			 		\vdots\\ 
			 		\displaystyle\frac{\partial f}{\partial \textbf{W}^{(L)}}^T 
				 \end{array} \right]
			 	 \left[ \frac{\partial f}{\partial \textbf{W}^{(1)}}, \ldots, \frac{\partial f}{\partial \textbf{W}^{(L)}} \right] = $$
			 	 
			 	 $$ = \left[ \begin{array}{ccc}
			 	 	\displaystyle\frac{\partial f}{\partial \textbf{W}^{(1)}}^T \frac{\partial f}{\partial \textbf{W}^{(1)}} & \ldots & \displaystyle\frac{\partial f}{\partial \textbf{W}^{(1)}}^T \frac{\partial f}{\partial \textbf{W}^{(L)}}\\
			  		\vdots & & \vdots\\
			  		\displaystyle\frac{\partial f}{\partial \textbf{W}^{(L)}}^T \frac{\partial f}{\partial \textbf{W}^{(1)}} & \ldots & \displaystyle\frac{\partial f}{\partial \textbf{W}^{(L)}}^T \frac{\partial f}{\partial \textbf{W}^{(L)}}\\
				 \end{array} \right],$$
where

				$$\frac{\partial f}{\partial \textbf{W}^{(\ell')}}^T \frac{\partial f}{\partial \textbf{W}^{(\ell'')}} = 
				\left[ \begin{array}{c}
					\displaystyle\frac{\partial f}{\partial \textbf{w}_1^{(\ell')}}^T\\
					\vdots\\ 
					\displaystyle\frac{\partial f}{\partial \textbf{w}_R^{(\ell')}}^T 
				\end{array} \right]
				\left[ \frac{\partial f}{\partial \textbf{w}_1^{(\ell'')}}, \ldots, \frac{\partial f}{\partial \textbf{w}_R^{(\ell'')}} \right] = $$
				
				$$ = \left[ \begin{array}{ccc}
					\displaystyle\frac{\partial f}{\partial \textbf{w}_1^{(\ell')}}^T \frac{\partial f}{\partial \textbf{w}_1^{(\ell'')}} & \ldots & \displaystyle\frac{\partial f}{\partial \textbf{w}_1^{(\ell')}}^T \frac{\partial f}{\partial \textbf{w}_R^{(\ell'')}}\\
					\vdots & & \vdots\\
					\displaystyle\frac{\partial f}{\partial \textbf{w}_R^{(\ell')}}^T \frac{\partial f}{\partial \textbf{w}_1^{(\ell'')}} & \ldots & \displaystyle\frac{\partial f}{\partial \textbf{w}_R^{(\ell')}}^T \frac{\partial f}{\partial \textbf{w}_R^{(\ell'')}}
				\end{array} \right].$$
			
				Let $\omega_{r' r''}^{(\ell)} = \langle \textbf{w}_{r'}^{(\ell)}, \textbf{w}_{r''}^{(\ell)} \rangle$ and assume, without loss of generality, that $1 \leq \ell' < \ell'' \leq L$. Thus we have that
			
				$$\frac{\partial f}{\partial \textbf{w}_{r'}^{(\ell')}}^T \frac{\partial f}{\partial \textbf{w}_{r''}^{(\ell'')}} = $$
			
				$$ = \left( \textbf{w}_{r'}^{(1)} \tilde{\otimes} \ldots \tilde{\otimes} \textbf{w}_{r'}^{(\ell'-1)} \tilde{\otimes} \textbf{I}_{I_{\ell'}} \tilde{\otimes} \textbf{w}_{r'}^{(\ell'+1)} \tilde{\otimes} \ldots \tilde{\otimes} \textbf{w}_{r'}^{(L)} \right)^T \cdot $$
				$$\left( \textbf{w}_{r''}^{(1)} \tilde{\otimes} \ldots \tilde{\otimes} \textbf{w}_{r''}^{(\ell''-1)} \tilde{\otimes} \textbf{I}_{I_{\ell''}} \tilde{\otimes} \textbf{w}_{r''}^{(\ell''+1)} \tilde{\otimes} \ldots \tilde{\otimes} \textbf{w}_{r''}^{(L)} \right) = $$
			
				$$ = \prod_{\ell \neq \ell', \ell''} \omega_{r' r''}^{(\ell)} \cdot \textbf{w}_{r''}^{(\ell')} \big( \textbf{w}_{r'}^{(\ell'')} \big)^T.$$
			
				In the case $\ell' = \ell''$ we have
			
				$$\frac{\partial f}{\partial \textbf{w}_{r'}^{(\ell')}}^T \frac{\partial f}{\partial \textbf{w}_{r''}^{(\ell')}} = $$
			
				$$ = \omega_{r' r''}^{(1)} \tilde{\otimes} \ldots \tilde{\otimes} \omega_{r' r''}^{(\ell'-1)} \tilde{\otimes} \  \textbf{I}_{I_{\ell'}}^2 \tilde{\otimes} \ \omega_{r' r''}^{(\ell'+1)} \tilde{\otimes} \ldots \tilde{\otimes} \omega_{r' r''}^{(L)} = $$ 
			
				$$ = \prod_{\ell \neq \ell'} \omega_{r' r''}^{(\ell)} \cdot \textbf{I}_{I_{\ell'}}.$$
			
				Therefore,
			
				$$\frac{\partial f}{\partial \textbf{W}^{(\ell')}}^T \frac{\partial f}{\partial \textbf{W}^{(\ell'')}} = 
				\left[ \begin{array}{ccc}
					\displaystyle \prod_{\ell \neq \ell', \ell''} \omega_{11}^{(\ell)} \cdot \textbf{w}_1^{(\ell')} \big( \textbf{w}_1^{(\ell'')} \big)^T & \ldots & \displaystyle \prod_{\ell \neq \ell', \ell''} \omega_{1R}^{(\ell)} \cdot \textbf{w}_R^{(\ell')} \big( \textbf{w}_1^{(\ell'')} \big)^T\\
					\vdots & & \vdots\\
					\displaystyle \prod_{\ell \neq \ell', \ell''} \omega_{R1}^{(\ell)} \cdot \textbf{w}_1^{(\ell')} \big( \textbf{w}_R^{(\ell'')} \big)^T & \ldots & \displaystyle \prod_{\ell \neq \ell', \ell''} \omega_{RR}^{(\ell)} \cdot \textbf{w}_R^{(\ell')} \big( \textbf{w}_R^{(\ell'')} \big)^T
				\end{array} \right]$$
when $\ell' \neq \ell''$. Finally, we have that

				$$\hspace{1.8cm} \frac{\partial f}{\partial \textbf{W}^{(\ell')}}^T \frac{\partial f}{\partial \textbf{W}^{(\ell')}} = 
				\left[ \begin{array}{ccc}
					\displaystyle \prod_{\ell \neq \ell'} \omega_{11}^{(\ell)} \cdot \textbf{I}_{I_{\ell'}} & \ldots & \displaystyle \prod_{\ell \neq \ell'} \omega_{1R}^{(\ell)} \cdot \textbf{I}_{I_{\ell'}}\\
					\vdots & & \vdots\\
					\displaystyle \prod_{\ell \neq \ell'} \omega_{R1}^{(\ell)} \cdot \textbf{I}_{I_{\ell'}} & \ldots & \displaystyle \prod_{\ell \neq \ell'} \omega_{RR}^{(\ell)} \cdot \textbf{I}_{I_{\ell'}}
				\end{array} \right]. \hspace{1.8cm}\square$$\bigskip

\section*{Acknowledgments}
I would like to acknowledge my advisor, professor Gregorio Malajovich, for the fruitful conversations we had and the help to put this paper in a reasonable format.

\end{document}


\maketitle

\section{A detailed example}

Here we include some equations and theorem-like environments to show
how these are labeled in a supplement and can be referenced from the
main text.
Consider the following equation:
\begin{equation}
  \label{eq:suppa}
  a^2 + b^2 = c^2.
\end{equation}
You can also reference equations such as \cref{eq:matrices,eq:bb} 
from the main article in this supplement.

\lipsum[100-101]

\begin{theorem}
  An example theorem.
\end{theorem}

\lipsum[102]
 
\begin{lemma}
  An example lemma.
\end{lemma}

\lipsum[103-105]

Here is an example citation: \cite{KoMa14}.

\section[Proof of Thm]{Proof of \cref{thm:bigthm}}
\label{sec:proof}

\lipsum[106-112]

\section{Additional experimental results}
\Cref{tab:foo} shows additional
supporting evidence. 

\begin{table}[htbp]
{\footnotesize
  \caption{Example table}  \label{tab:foo}
\begin{center}
  \begin{tabular}{|c|c|c|} \hline
   Species & \bf Mean & \bf Std.~Dev. \\ \hline
    1 & 3.4 & 1.2 \\
    2 & 5.4 & 0.6 \\ \hline
  \end{tabular}
\end{center}
}
\end{table}

\bibliographystyle{siamplain}
\bibliography{references}